\documentclass[12pt]{amsart}

\usepackage{amssymb}

\usepackage{enumerate}

\usepackage{graphicx}

\makeatletter
\@namedef{subjclassname@2010}{%
  \textup{2010} Mathematics Subject Classification}
\makeatother

\newtheorem{thm}{Theorem}[section]
\newtheorem{cor}[thm]{Corollary}
\newtheorem{lem}[thm]{Lemma}
\newtheorem{prop}[thm]{Proposition}

\theoremstyle{definition}
\newtheorem{defin}[thm]{Definition}
\newtheorem{rem}[thm]{Remark}

\newtheorem{algo}[thm]{Algorithm}
\newtheorem{claim}[thm]{Claim}

\numberwithin{equation}{section}

\newcommand{\lineclear}
       {\rule{0pt}{0pt}\nopagebreak\par\nopagebreak\noindent}

\def\hide#1{}

\newcommand{\N}{\mathbb N}
\newcommand{\Nplus}{{\mathbb N^*}}

\newcommand{\Z}{\mathbb {Z}}
\newcommand{\R}{\mathbb {R}}
\newcommand{\C}{\mathbb {C}}
\newcommand{\D}{\mathbb D}
\newcommand{\J}{\mathcal J}
\newcommand{\K}{\mathcal K}
\newcommand\I{{\mathcal I}}

\newcommand{\Cbar}{\overline{\C}}
\newcommand\M{\mathcal M}
\newcommand{\Mandel}{\M}

\newcommand{\Circle}{{\mathbb S}^1}
\newcommand{\ie}{{\it i.e.,\ }}
\newcommand{\eg}{{\it e.g.\ }}

\newcommand\sm{\setminus}
\newcommand\ovl[1]{\overline {#1}}
\newcommand{\0}{{\tt 0}}
\newcommand{\1}{{\tt 1}}
\newcommand{\2}{{\tt 2}}
\renewcommand{\*}{{\tt \star}}
\newcommand{\Sym}{\Sigma^1}
\newcommand{\Symg}{\Sigma}
\newcommand{\Syms}{\Sigma^\star}

\newcommand{\IntAddr}{\to}
\newcommand{\IntAdr}{\to}
\newcommand{\makehigh}{\rule{0pt}{12pt}}
\newcommand{\orb}{\mbox{\rm orb}}
\newcommand{\diam}{\mbox{\rm diam}}
\newcommand{\hdim}{\mbox{\rm dim}_H}
\newcommand{\rhox}{\rho_{\nu, x}}

\newcommand{\Biac}{{\mathcal Biac}}
\newcommand{\evil}{evil}
\newcommand{\tame}{tame}

\renewcommand\theta{\vartheta}
\newcommand{\eps}{\varepsilon}
\newcommand{\z}{\zeta}

\newcommand{\step}{\mbox{{\sc Step}}}  
\newcommand{\binary}{b}
\newcommand{\doub}{D}
\renewcommand{\phi}{\varphi}

\newcommand{\pair}[2]{\{#1,#2\}}
\renewcommand{\pair}[2]{\langle #1,#2\rangle }

\newcommand{\Biack}{E_k}
\newcommand{\heading}[1]{\par\medskip\noindent\textbf{#1}}

\newcommand{\diskbar}{\ovl{\D}}

\begin{document}
\baselineskip=17pt 

\title[Hausdorff dimension of biaccessible angles \today]
{Hausdorff dimension of biaccessible angles \\[1mm]
for quadratic polynomials.}
\author{Henk Bruin and Dierk Schleicher}

\address{
Faculty of Mathematics,
University of Vienna,
Oskar Morgensternplatz 1,
1090 Vienna, Austria}
\email{henk.bruin@univie.ac.at}

\address{Jacobs University Bremen, Research I,
P.O. Box 750 561, D-28725 Bremen,
Germany}
\email{dierk@jacobs-university.de}



\begin{abstract}
A point $c$ in the Mandelbrot set is called biaccessible if 
two parameter rays land at $c$. 
Similarly, a point $x$ in the Julia set of a polynomial 
$z \mapsto z^2+c$ is called biaccessible if two dynamic rays land at $x$.
In both cases, we say that the external angles of these two rays 
are biaccessible as well.

In this paper we describe a purely combinatorial characterization of 
biaccessible (both dynamic and parameter) angles,
and use it to give detailed estimates of the Hausdorff dimension
of the set of biaccessible angles.
\end{abstract}

\maketitle

\section{Introduction}\label{sec:intro}

Dynamic rays and their landing properties are a key tool to
understanding (the topology of) Julia sets of polynomials. In particular, the structure of the Julia set is determined by rays that land at a common point: at least in good cases (under the assumption of local connectivity), the knowledge of which rays land together gives a homeomorphic model for the Julia set that is known as Douady's \emph{pinched disk model} \cite{DouadyCompacts}. Very similarly, Thurston developed his concept of \emph{invariant laminations} \cite{ThurstonLaminations} that provides a uniform topological model at least of quadratic polynomials, based upon a single quantity, the external angle, which determines which rays land together (at least combinatorially). Analogous statements hold for the Mandelbrot set, the parameter space of quadratic polynomials with connected Julia sets: there is a simple topological model, the \emph{quadratic minor lamination}, that describes the Mandelbrot set in terms of which rays (should) land together. In all these cases, a point in $\C$ is called \emph{biaccessible} if it is the landing point of two or more rays. 

Biaccessibility is of interest from several more points of view.
For instance, mating constructions \cite{rees, shishikura, tanlei2}
and certain constructions of space-filling curves \cite{sirvent1, sirvent2} 
also rely on the biaccessibility of external angles
(and rays).  

It is of interest to quantify ``how many'' rays land together. The \emph{biaccessibility dimension} is defined as the Hausdorff dimension of those external angles for which the corresponding rays land together with some other ray. Thurston observed \cite{thurston} that at least for postcritically finite parameters the biaccessibility dimension equals (up to a factor of $\log 2$) the \emph{core entropy} of the polynomial, \ie the topological entropy of the restriction to the Hubbard tree. 
This relation holds true in greater generality, for instance for postcritically infinite parameters with a compact Hubbard tree \cite{tiozzo}. For an appropriately extended definition of core entropy for arbitrary quadratic polynomials with connected Julia set, this relation is in fact true in full generality; see Jung's appendix in \cite{DS}. Inspired by questions of Thurston and Hubbard, the question of continuity of core entropy is a topic of core interest: see \cite{thurston,CT,tiozzo,jung,TanLeiEntropy} and the recent papers \cite{tiozzo2,DS}.

Biaccessible angles have been studied in terms of Lebesgue measure on 
$\Circle$, in particular by
Smirnov \cite{smirnov1} and Zdunik \cite{zdunik}. For any polynomial 
Julia set that is not an interval, the set of biaccessible angles on 
$\Circle$ has $1$-dimensional Lebesgue measure zero. In other words, 
the biaccessible points have harmonic measure zero. 
This was strengthened in  \cite{MeerkampSchleicher} proving that the biaccessibility dimension (the Hausdorff dimension of the set of biaccessible dynamic angles) is strictly less than $1$ (except, of course, when the Julia set is an interval).
For further results on different aspects of biaccessibility, see for 
instance Zakeri \cite{zakeri3} and \cite{SchleicherZakeri}.

In this paper, we will take a purely combinatorial point of view.
That is, we express {\em combinatorial biaccessibility}
in terms of external angles, the angle doubling map, and their itinerary
with respect to the partition of $\Circle$ defined by $\theta/2$ and 
$(1+\theta)/2$, where $\theta$ is the parameter angle. 
We denote by $\Biac_\theta$ the set of combinatorially biaccessible 
$\varphi \in \Circle$, and since this condition carries over to parameter 
space, we can define $\Biac$ analogously as the set of 
combinatorially biaccessible parameter angles $\theta \in \Circle$.
To summarize our main results, we 
\begin{itemize}
\item give detailed estimates of
$\dim_H(\Biac_\theta)$ (Section~\ref{sec:comb_dim_sequences}) 
and $\dim_H(\Biac)$ (Section~\ref{Sec:DimParaAngles});
\item treat the case of real parameters $c \in \Mandel \cap \R$,
and conclude that $\dim_H(\Biac) = 1$, which is due to the angle $\theta = \frac12$: away from any neighborhood of $\theta = \frac12$, the dimension
is strictly smaller than $1$ 
(see the remark below Theorem~\ref{ThmHdim_mandel2});
\item 
describe exactly those angles $\theta$ for which $\dim_H(\Biac_\theta)=0$. 
\end{itemize}
We give precise statements of these results in Section~\ref{sec:main}. The necessary combinatorial language will be developed in Section~\ref{sec:combinatorial}; in particular, we relate topological and combinatorial biaccessibility.

We emphasize that the Hausdorff dimension estimates 
obtained in this paper are for sets of external angles, \ie 
subsets of $\Circle$.
As far as we are aware, there is no direct relation to the dimension
of the set of biaccessible points either in $\J_c$ or in $\Mandel$.
For instance, Lyubich \cite{Lyubich} showed that the
parameters in $\Mandel \cap \R$ representing infinitely renormalizable maps
form a set of Hausdorff dimension at least $ \frac12$.
In contrast, 
the Hausdorff dimension of infinitely renormalizable parameter angles
is zero (see Section~\ref{Sub:RenormalizableAngles}).

\looseness-1
Taking a purely combinatorial approach bypasses the complications
of non-locally connected Julia sets (and potentially the Mandelbrot set).
Julia sets and the Mandelbrot set are possibly not locally connected, in which 
case some external rays may not land, or not land where they are 
``combinatorially'' supposed to land. We will show in Proposition~\ref{Prop:ExtraAngles} that this has no impact on the Hausdorff dimension of the angles of biaccessible points.

\medskip
\emph{Acknowledgement}. The authors thank Wolf Jung and Marten Fels 
for their critical reading of this text and related manuscripts,
and the referee(s) for their insightful comments.
Also, we would like to thank the Erwin-Schr\"odinger-Institut in Vienna 
(specifically, the workshop ``Ergodic Theory and Holomorphic Dynamics''
September--October 2015) for their support during an important phase of this work.

\section{Statements of Main Results}
\label{sec:main}

\subsection{Rays and biaccessibility}

\looseness-1
We start by giving a quick review of rays and their landing properties, before we relate this, in later subsections, to combinatorial questions. 
We will only consider quadratic polynomials 
$p_c(z) = z^2+c$ on the Riemann sphere $\Cbar$.
In this setting, $\infty$ is a super-attractive fixed point.
All points that converge to $\infty$ under iteration of $p_c$
belong to the basin $A(\infty)$; the remaining points belong to
the {\em filled-in Julia set} $\K_c=\C\sm A(\infty)$. The {\em Julia set}
$\J_c$ is the common boundary of $\K_c$ and $A(\infty)$.

If $\K_c$ is connected, then there is a unique Riemann map $\psi_c\colon\Cbar\sm\K_c\to\Cbar\sm\diskbar$ with 
$\psi_c(\infty) = \infty$ and $\psi_c(z)/z\to 1$ as $z\to\infty$.  \emph{B\"ottcher coordinates} on $A(\infty)$ are defined as preimages of polar coordinates on $\Cbar\sm\diskbar$: every $z\in\C\sm\K_c$ has its \emph{potential} $|\psi_c(z)|$ and its \emph{external angle} $\arg\psi_c(z)/2\pi$ (so that external angles are measured in terms of full turns, where the full circle has measure $1$, not $2\pi$ in radians). 

The B\"ottcher map $\psi_c$ conjugates
$p_c$ to the map $z \mapsto z^2$ as follows:  $(\psi_c(z))^2 = \psi_c \circ p_c(z)$ for all $z \in A(\infty)$.
Given an angle $\varphi \in [0,1)$, the {\em dynamic ray at angle $\varphi$} is
the set $R_c(\varphi) := \psi_c^{-1}(\{ re^{2\pi i \varphi} : r > 1\})$, and the ray is said to {\em land}
if $\lim_{r \to 1}\psi_c^{-1}(re^{2\pi i \varphi})$ exists; this limit
is called the {\em landing point} (it is always in $\J_c$).
An external angle is called {\em biaccessible} if there is another external angle so that the associated rays have the same landing 
points.

The situation in parameter space is analogous. 
The Mandelbrot set $\M$ is defined as the set of parameters $c$ for which the filled-in Julia set $\K_c$ is connected. 
There is a Riemann map $\psi:\Cbar \sm \Mandel \to \Cbar \sm \diskbar$
for the exterior of the Mandelbrot 
set; in these terms, we can define parameter rays 
$R(\theta) := \psi^{-1}(\{ re^{2\pi i \theta} : r > 1\})$
and study their landing points and
biaccessibility. The Riemann map $\psi$ and parameter rays were introduced by Douady and Hubbard in their Orsay Notes 
\cite{Orsay}.

\subsection{Itineraries, sequences and the $\rho$-function}
\begin{figure}[htbp]
\includegraphics{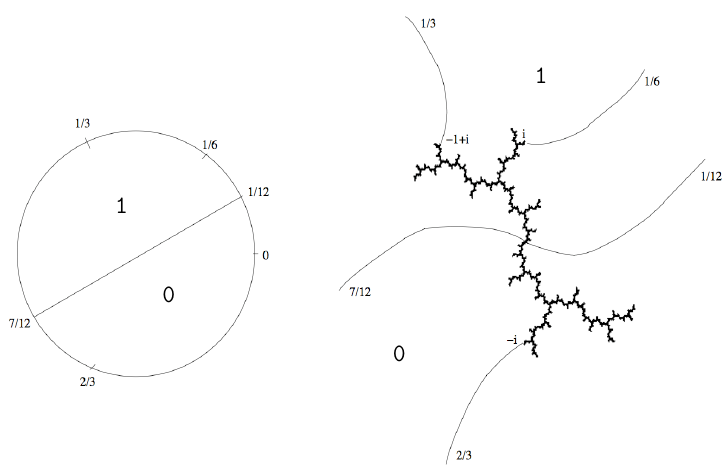}
\caption{Left: the kneading sequence of an external angle $\theta$
(here $\theta = 1/6$)
is defined as the itinerary of the orbit of $\theta$ under angle
doubling, taken with respect to the
partition formed by the angles $\theta/2$ and $(\theta+1)/2$.
Right: in the dynamics of a polynomial for which the $\theta$-ray
lands at the critical value, an analogous partition is formed by
the dynamic rays at angles $\theta/2$ and $(\theta+1)/2$, which
land together at the critical point. }
\label{FigDefKneading}
\end{figure}

\looseness-1
Our starting point is an angle $\theta \in \Circle = \R / \Z$ that we view as external parameter. It is used to partition $\Circle$ and define symbolic dynamics
for the angle doubling map $\doub:\Circle \to \Circle$, $\varphi \mapsto 2\varphi \pmod 1$.
\begin{defin}{(Itinerary and Kneading Sequence of External Angle).}
\label{DefKneading} \lineclear
Given an external angle $\theta\in\Circle$, we associate
to each $\varphi\in\Circle$ its {\em itinerary}
$\nu_\theta(\varphi)= \nu_1\nu_2\ldots$ with $\nu_k\in\{\0,\1,\*\}$ by:
\[
\nu_k:=\left\{
\begin{array}{ll}
\0 \qquad & \mbox{if }\ \doub^{\circ k-1}(\varphi) \in \left(\frac{1+\theta}{2} , \frac{\theta}{2}\right), \\[1mm]
\1 & \mbox{if }\ \doub^{\circ k-1}(\varphi) \in \left(\frac{\theta}{2}, \frac{1+\theta}{2} \right), \\[1mm]
\* & \mbox{if }\ \doub^{\circ k-1}(\varphi)\in\left\{\frac{\theta}{2},\frac{1+\theta}{2} \right\},
\end{array}
\right.
\]
where the intervals are interpreted with respect to cyclic
order. The {\em kneading sequence}
$\nu(\theta)$ of $\theta$ is
its itinerary with respect to itself:
$\nu(\theta):=\nu_\theta(\theta)$; see Figure~\ref{FigDefKneading}.

Finally, for $\phi\in\Circle$ we say that $\step(\phi)=k$ if $D^{\circ k}(
\phi)=\theta$ and $k\ge 0$ is minimal with this property. 
\end{defin}

Note that in the definition of the kneading sequence, a change in
$\theta$ amounts to a change in the orbit $\orb_\doub(\theta)$,
as well as a change in the partition itself. Our definition involves the convention that every kneading sequence $\nu(\theta)$ starts with the symbol $\1$, except for $\theta=0$.

A kneading sequence $\nu$ contains a $\*$ at position $n$
if and only if $\theta$ is periodic with period $n$; the exact period 
of $\theta$ may divide $n$. (There are also non-periodic angles that yield periodic 
kneading sequences without $\*$.) We say that a
sequence $\nu$ is {\em $\*$-periodic of period $n$} if $\nu =
\ovl{\nu_1 \ldots \nu_{n-1}\*}$ with $\nu_1=\1$ and
$\nu_i \in \{ \0,\1 \}$ for $1 < i < n$; this happens if and only if $\theta$ is periodic of exact period $n$.
Write $\N = \{ 0, 1, 2, 3, \dots\}$ and $\Nplus = \{ 1, 2, 3, \dots\}$.
Let
\begin{eqnarray*}
\Symg &:=& \{\0,\1\}^\Nplus \,\,, \\
\Sym &:=&
\{\nu\in\Symg\colon \mbox{the first entry in $\nu$ is $\1$}\}  \,\,,\\
\Syms &:=& \Sym\cup \{\mbox{all $\*$-periodic sequences except 
$\ovl\*$}\} \,\,.
\end{eqnarray*}
The metric on $\Symg$ is $d(x,y) = \sum_{i=1}^\infty 2^{-i} |x_i-y_i|$.
In order to avoid silly counterexamples, $\ovl{\*}$ is not
considered to belong to $\Syms$. All sequences in $\Syms$
will be called {\em kneading sequences}, regardless of
whether or not they occur
as the image of an angle $\theta\in\Circle$. 

To compress the information of the kneading sequence $\nu$,
it is useful to introduce the $\rho$-function.

\begin{defin}{($\rho$-Function and Internal Address\footnote{
The $\rho$-function is of fundamental importance in the work of
Penrose~\cite{Pe} under the name of {\em non-periodicity function};
the internal address is called {\em principal non-periodicity function}.}).}
\label{DefRho} \lineclear
For a sequence $\nu \in\Syms$, define
\[
\rho_{\nu}:\Nplus \to \Nplus\cup\{\infty\}, \quad
\rho_\nu(n) = \inf \{ k > n: \nu_k \neq \nu_{k-n} \}.
\]
We usually write $\rho$ for $\rho_\nu$ and call
$\orb_\rho(k) = \{ \rho^{\circ i}(k) \}_{i \geq 0}$
the $\rho$-orbit of $k$. The case $k=1$ is the most important one;
this is the {\em internal address} of $\nu$
and we denote it as
\[
1 =  S_0 \IntAddr S_1 \IntAddr S_2  \IntAddr \dots
\]
\end{defin}
\noindent
The name internal address is motivated by the fact that
if $\theta$ is the parameter angle of $c \in \partial \Mandel$ and
$\gamma \subset \Mandel$ is a (combinatorial) arc connecting $0$ to $c$,
then the periods of the hyperbolic components of $\Mandel$ intersecting $\gamma$ having
a lower period than any later hyperbolic component intersecting $\gamma$
form precisely the entries of the internal address, see \cite{IntAdr,IntAddr}.

The map from kneading sequences in $\Sym$ to internal addresses is
injective. In fact, the algorithm of this map can easily be inverted:

\begin{algo}{(From Internal Address to Kneading Sequence).}
\label{AlgIntAdrKneading} \lineclear
The following inductive algorithm turns internal addresses into
kneading sequences in $\Sym$: the internal address $S_0=1$ has
kneading sequence $\ovl\1$, and given the kneading sequence $\nu^k$
associated to $1\IntAdr S_1\IntAdr\ldots\IntAdr S_k$,
the kneading sequence associated to $1\IntAdr S_1\IntAdr
\ldots\IntAdr S_k\IntAdr S_{k+1}$ consists of the first
$S_{k+1}-1$ entries of $\nu^k$, followed by the opposite to the entry
$S_{k+1}$ in $\nu$ (switching $\0$ and $\1$), and then repeating
these $S_{k+1}$ entries periodically.
\end{algo}

\begin{proof}
The kneading sequence $\ovl \1$ has internal address $1$. If $\nu^k$
has internal address $1\IntAdr S_1\ldots\IntAdr
S_k$ and $\nu$ is the kneading sequence of period $S_{k+1}$ as
constructed in the algorithm, then the internal address of $\nu$
clearly starts with $1\IntAdr S_1\IntAdr\ldots\IntAdr
S_k$, and $\rho_{\nu}(S_k)=S_{k+1}$, so the internal address of
$\nu$ is $1\IntAdr S_1\IntAdr \ldots \IntAdr S_k\IntAdr
S_{k+1}$.
\end{proof}

\subsection{Hausdorff Dimension of Biaccessible Angles}

Our results are about biaccessible points both in Julia sets and in the Mandelbrot set. 
If $K\subset\C$ is compact, then a point $z\in\partial K$ is \emph{topologically biaccessible} 
if there are two curves $\gamma_1,\gamma_2\colon[0,1]\to\Cbar$ with $\gamma_i(0)=\infty$, $\gamma_i(1)=z$, $\gamma_i([0,1))\subset\Cbar\sm K$ (for $i=1,2$) and so that $\gamma_1$ and $\gamma_2$ are not homotopic in $\Cbar\sm K$ fixing endpoints. We will translate this into a combinatorial setting in Section~\ref{sec:combinatorial} but first state our main results here. These results will be given in terms of two quantities $N$ and $\kappa$ that we define first.

(i) For a kneading sequence $\nu = \1\nu_2\nu_3\nu_4\dots \in \Sym$, let
\begin{equation}\label{EqN_of_nu}
N = N(\nu) := 1 + \min \{ i > 1 \colon \nu_i = \1\},
\end{equation}
so that $N(\nu)-1$ is the position of the second $\1$ in $\nu$ (and 
$N(\nu)=\infty$ if $\nu=\1\ovl{\0}$); hence $N(\nu)\ge 3$.
Set $M:=\lfloor N/2\rfloor-1$ and
\begin{align}\label{Eq_U1N}
U_1(N)  &:= 
1-\frac{1}{2^N\log 2^N} \qquad \text{and}
\\
L_1(N) &:= \left\{
\begin{array}{ll}
1-\frac{1}{2^M \log 2^{M} }
  & \text{ if } N\ge 6\\
1/2  & \text{ if } N = 5, \\
0 & \text{ if } N\in\{3,4\}.
\end{array} \right.
\end{align}
Finally, set $L_1(N)=U_1(N)=1$ if $N=\infty$.

(ii)  For a kneading sequence $\nu$ with internal address
$1 \IntAddr S_1\IntAddr S_2 \IntAddr \dots$, let
\begin{equation}\label{Eqkappa_of_nu}
\kappa = \kappa(\nu) := \sup\{ k \geq 1 : S_j \text{ is a multiple of } S_{j-1}
\text{ for all } 1 \leq j \leq k \} \;.
\end{equation}
We define $L_2(S_\kappa)$ and $U_2(S_\kappa)$ as follows:

\begin{itemize}
\item[(a)]
if $\kappa=\infty$, then $L_2(S_\kappa)=U_2(S_\kappa)=0$;
\item[(b)]
if $\kappa<\infty$ and $\nu$ is periodic of period $S_\kappa$, then again $L_2(S_\kappa)=U_2(S_\kappa)=0$;
\item[(c)] \label{caseC}
otherwise (\ie if $\kappa<\infty$ and $\nu$ is not periodic of period $S_\kappa$, hence $S_{\kappa+1}<\infty$), we set 
\begin{equation}\label{Eq_U2kappa}
L_2(S_\kappa) := {1}/{S_{\kappa+1}}, \quad
U_2(S_\kappa) := \sqrt{{7}/{S_{\kappa+1}}}
\;.
\end{equation}

\end{itemize}
With these definitions, define the interval
\begin{equation}
\I(N,S):= \left[\ \max\left\{L_1(N), L_2(S)\rule{0pt}{10pt}\right\} \ , \
\min\left\{U_1(N), U_2(S)\rule{0pt}{10pt}\right\}\ \makehigh \right].
\label{EqDimInterval}
\end{equation}

We define combinatorial biaccessibility in Definition~\ref{Def:CombBiac}; this leads to the following two sets (for the dynamical planes respectively parameter space) that we will investigate:
\[
\Biac_\theta = \{ \varphi \in \Circle \colon \text{$\phi$ is combinatorially biaccessible with respect to $\theta$} \},
\]
and
\[
\Biac = \{ \theta \in \Circle \colon \text{$\theta$ is combinatorially biaccessible} \}.
\]
With these definitions, our first main result is as follows.

\begin{thm}{(Hausdorff Dimension of Biaccessible Angles).}
\label{ThmHdimBiaccAngles} \lineclear
For every parameter angle $\theta\in\Circle$, 
$$
\dim_H(\Biac_\theta) \in \I(N,S)
$$ 
for $N = N(\nu(\theta))$ and $S = S_\kappa(\nu(\theta))$. 
In particular, the set of biaccessible external dynamic angles
has Hausdorff dimension less than $1$ unless $\theta=1/2$.
\end{thm}

This implies that the harmonic measure of the biaccessible
points in quadratic Julia sets is zero, unless $\theta=\frac12$; 
this was of course known earlier \cite{smirnov1,zakeri3,zdunik}. Our result that the Hausdorff dimension of biaccessible angles is less than $1$ except when the Julia set is an interval (which, in the quadratic case, means $\theta=1/2$) was later generalized to all degrees in \cite{MeerkampSchleicher}.

The similarity of Theorem~\ref{ThmHdimBiaccAngles} and the next one
underlines the similarity of the structure of $\J_c$ and the 
local structure of $\Mandel$ near $c$ as in \cite{tanlei1}.

\begin{thm}{(Hausdorff Dimension of Biaccessible Parameter Angles).}
\label{ThmHdim_mandel2}\lineclear
For any $N \ge 3, \kappa \ge 1$ we have
\[
\dim_H\left(\Biac \cap \left\{ \theta \in \Circle : N(\nu(\theta)) = N \text{ and } 
S_\kappa(\nu(\theta)) = S \right\}\right) \in \I(N,\kappa).
\]
\end{thm}

\begin{rem}
The estimate $[L_2(S_\kappa),U_2(S_\kappa)]$ is best where the biaccessibility dimension is small. 
In every case we either have $L_2(S_\kappa)=U_2(S_\kappa)=0$ or $U_2(S_\kappa)>L_2(S_\kappa)>0$, 
so that we specify exactly in which situation the biaccessibility dimension is zero: 
this happens if and only if the parameter is on the closed main molecule 
(see Proposition~\ref{PropDimBiacRenorm}). 
Our result implies continuity of the biaccessibility dimension on the closed main molecule of $\M$: if $(c_n)$ is a sequence of parameters in $\M$ that converges to the main molecule and with associated external angles $\theta_n$, then $\dim_H(\Biac_{\theta_n})\to 0$.

Similarly, the estimate $[L_1(N),U_1(N)]$ is especially good near parameters where the biaccessibility dimension is maximal, that is near the ``antenna tip'' at $c=-2$. Again, we have $L_1(N)=U_1(N)=1$ or $L_1(N)<U_1(N)<1$, so we specify exactly in which situation the dimension is $1$: this happens if and only if $c=-2$, which was generalized in \cite{MeerkampSchleicher}. Again our result implies continuity in the approach to this point.
\end{rem}

\begin{rem}
Since, as mentioned before, core entropy (appropriately defined for all polynomials in $\M$) 
equals biaccessibility dimension times by $\log 2$, these results give a different proof 
for continuity of the core entropy at the main molecule and at $c=-2$ (\ie for parameters with 
core entropy in $\{0,\log 2\}$) with precise estimates. (Of course, continuity everywhere 
was proved in \cite{tiozzo2,DS}.) 
\end{rem}

\begin{rem} In particular, 
the set of biaccessible parameter angles has Hausdorff dimension $1$ 
but Lebesgue measure zero: outside of every neighborhood of $\frac12$ they have  Hausdorff dimension less than $1$.
The same holds for the set of parameter angles with landing point
on the real antenna $\Mandel \cap \R$.
This follows because the collection
of kneading sequences $\nu$ of the form used in our proof (formula
\eqref{eq:BW} to be precise)
has the property that $\orb_\rho(1) \cap \orb_{\rho}(N-1) = \emptyset$.
According to \cite{thunberg}, this means that $\nu$ is the kneading sequence
of a real quadratic map, so the part of the proof below that refers to
\eqref{eq:BW} (\ie biaccessible parameter angles close to $\frac12$)
automatically gives that set of ``real'' parameter angles
has indeed Hausdorff dimension $1$. We would also like to mention recent work by Tiozzo especially on the kneading sequences and external angles on $\M\cap\R$ \cite{tiozzo}: for every $c\in\M\cap\R$, the set of external angles of $\M$ that lands on $[c,0]\subset(\M\cap\R)$ 
has the same Hausdorff dimension as the set of all biaccessible angles of $\J_c$.
\end{rem}

\subsection{Renormalizable Angles}
\label{Sub:RenormalizableAngles}
Now we draw some direct consequences of the methods of the main results that have to do with renormalization.

\begin{defin}{(Renormalizable).}
\label{DefRenormalizable} \lineclear
A quadratic polynomial $p_c:z \mapsto z^2+c$ is {\em $M$-renormalizable}
for $M \ge 2$ if there exist neighborhoods $U \subset V$ of the critical point $0$, with $U$ compactly contained in $V$, 
such that $p_c^{\circ M}:U \to V$ is a degree $2$ branched covering such that
$p_c^{\circ iM}(0) \in U$ for all $i \geq 1$. The integer $M$ is called the {\em period of renormalization}. 
The set $K_U:=\{z\in U\colon p_c^{\circ M i}(z)\in U \text{ for all } 
i\geq 0\}$ is called the \emph{little filled-in Julia set} of the 
renormalization.
The renormalization is called \emph{simple} if $K_U\cap p_c^{\circ 
i}(K_U)$ does not disconnect $K_U$ for $i=1,\dots,M-1$; otherwise it is called a {\em crossed renormalization}. 
\end{defin}

If $p_c$ is simple renormalizable of period $M$, then $M = S_k$ is 
an entry in the internal address and all successive entries
$S_j$, $j \geq k$, are multiples of $M$, and conversely 
\cite{IntAdr,IntAddr}.
The corresponding kneading sequence has the form
\begin{equation}\label{eq:nurenorm}
\nu = \nu(\theta) = \nu_1 \nu_2 \dots \nu_{M-1} \nu_M
 \nu_1 \nu_2 \dots \nu_{M-1} \nu_{2M}  \nu_1 \nu_2 \dots \nu_{M-1} \nu_{3M} \dots
\end{equation}
where either $\nu_M\nu_{2M}\nu_{3M}\dots$ or its opposite sequence 
 $\nu'_M\nu'_{2M}\nu'_{3M}\dots$ (where
$\nu'_i = \1-\nu_i$) is the kneading sequence
of the renormalization $p_c^{\circ M}$.

We recover a result by Manning on external angles \cite[page 523]{Manning} and extend it to kneading sequences as follows.

\begin{prop}{(Dimension of Renormalizable Angles and Kneading Sequences).}
\label{PropDimRenorm} \lineclear
For any $M \geq 2$, the Hausdorff dimension of the set of 
simple $M$-renormalizable parameter angles as well as 
of the set of simple $M$-renormalizable 
kneading sequences is at most $1/M$. 
The Hausdorff dimension of infinitely renormalizable parameter angles 
as well as of infinitely renormalizable kneading sequences is $0$.
\end{prop}

Let some angle $\theta$ have internal address
$1 \IntAddr S_1 \IntAddr S_2 \IntAddr \dots$ or
$1 \IntAddr S_1 \IntAddr S_2 \IntAddr \dots \IntAddr S_k$ (if it is finite).
If $S_j$ is a multiple of $S_{j-1}$ for all $1 \le j < \infty$
(or  $1 \le j \le k$),
then we say that $\theta$ is associated to the \emph{main molecule of the Mandelbrot set}. 
The angle may be non-renormalizable (if the internal address is just the single entry $1$), finitely renormalizable (if the internal address has finitely many entries, but at least $2$), or infinitely renormalizable (if the internal address is infinite). 

The best known example of an infinitely renormalizable map from the main molecule
is the Feigenbaum-Coullet-Tresser map $p_c$ with $c = -1.4011551890...$,
where $S_j = 2^j$.

\begin{prop}{(Main Molecule Angles).}
\label{PropDimBiacRenorm} \lineclear
The following two conditions are equivalent for a parameter angle $\theta$:
\begin{itemize}
\item $\dim_H(\Biac_\theta)=0$;
\item $\theta$ is associated to the main molecule of $\M$.
\end{itemize}
\end{prop}

\section{The Combinatorial Approach}\label{sec:combinatorial}

\subsection{The Hubbard Tree and non-admissible kneading sequences}

The Hubbard tree of a postcritically finite polynomial is defined as the connected hull of 
the union of all critical orbits within the filled-in Julia set (subject to a regularity 
condition on how to pass through bounded Fatou components). We view a Hubbard tree 
of a quadratic polynomial as a finite abstract tree $T$ with dynamics $f\colon T\to T$ subject to the following conditions:
\begin{enumerate}
\item $f:T \to T$ is a local embedding, except at a single critical point $c_0$.
\item This critical point divides $T$ into (at most) two parts, labeled $\1$ (the part containing the critical value $c_1 = f(c_0)$) and $\0$, while $c_0$ itself gets the symbol $\star$.
Using these symbols, we can define  itineraries in the usual way.
\item The endpoints of $T$ lie on the critical orbit.
\item All {\em marked points} (\ie branch points and points on the critical orbit) have distinct itineraries (with respect to the partition introduced by $c_0$).
\end{enumerate}
It was shown in \cite{BKS} that for every $\*$-periodic or preperiodic $\nu$,
there is a Hubbard tree $(T,f)$ such that $\nu$ is the itinerary of the critical value.
Moreover, $\nu$ uniquely determines the dynamics on the marked points,
and their numbers of arms, so the pair $(T,f)$ is determined uniquely as an abstract tree with dynamics (up to homotopy relative to the marked points).
However, there can be multiple ways (up to homotopy) of embedding $T$ into $\C$
such that the dynamics $f$ extends to local homeomorphism on $\C$.
This depends on how the arms of a branch point $p$ are arranged;
if $p$ is periodic of period $n$, then $f^{\circ n}$ permutes the arms of 
$p$ in a transitive way, but the choice of \emph{combinatorial rotation number} is still free. This is the reason why multiple parameter angles lead to the same kneading sequence (see also Figure~\ref{FigNu}) and the same internal address. In order to distinguish different embeddings of these trees, \emph{angled internal addresses} are required \cite{IntAddr}. We will come back
to this in Section~\ref{Sec:DimParaAngles}, where we estimate the total
number of different embeddings in $\C$.

Not every sequence in $\Sym$ or $\Syms$ occurs as the kneading sequence of an 
external angle, or of a quadratic polynomial, not even every periodic sequence. The main result of 
\cite{BS1} is an explicit condition that states which sequences occur: a kneading 
sequence or internal address does occur unless it fails the following 
admissibility condition for some period $m$.
Failing the condition forces the Hubbard tree to have a periodic branch point of period $m$ with certain specific properties; such periodic orbits are called \emph{evil}.

\begin{defin}{(The Admissibility Condition).}
\label{DefAdmissCond} \lineclear
A kneading sequence $\nu\in\Syms$ {\em fails the admissibility
condition for period $m$} (the {\em evil} period) if the following three conditions hold:
\begin{enumerate}
\item
the internal address of $\nu$ does not contain $m$;
\item
if $k<m$ divides $m$, then $\rho(k)\leq m$;
\item
$\rho(m)<\infty$ and if $r\in\{1,\ldots,m\}$ is congruent to 
$\rho(m)$ modulo $m$, then
$\orb_\rho(r)$ contains $m$.
\end{enumerate}
A kneading sequence {\em fails the admissibility condition} if it 
does so for some $m \geq 1$.

\noindent
An internal address {\em fails the admissibility condition} if its
associated kneading sequence (from Definition~\ref{DefRho} and Algorithm~\ref{AlgIntAdrKneading}) does.
\end{defin}

\begin{rem}
The main result in \cite{BKS} is that for every $\*$-periodic kneading sequence $\nu$ there exists a unique Hubbard tree (without embedding into $\C$) for which the critical value has kneading sequence $\nu$ (with respect to the unique critical point of this tree), and the main result of \cite{BS1} is that this tree can be embedded into the plane so that the dynamics is compatible with this embedding if and only if $\nu$ 
does not fail this admissibility condition. Since every tree thus embedded is realized by a quadratic polynomial and has two characteristic periodic dynamic rays, this means that a $\*$-periodic kneading sequence is realized by an external angle if and only if the kneading sequence does not fail the admissibility condition.
\end{rem}

\subsection{Combinatorial biaccessibility} 
In this section we give the central
combinatorial characterization of biaccessibility. We denote an ordered pair of angles by $\pair{\phi}{\phi'}$ (this is technically the same as the more usual notation $(\phi,\phi')$ for an element of $\Circle\times\Circle$, but we will reserve the latter for the open interval). We say that an angle pair $\pair{\phi}{\phi'}$ \emph{separates} another angle pair $\pair{\tilde\phi}{\tilde\phi'}$ if $\tilde\phi$ and $\tilde\phi'$ are in different components of $\Circle\sm\{\phi,\phi'\}$. 

Suppose that in some dynamical plane the dynamic ray at angle $\theta$ lands at the critical value (a discussion that holds in more general situations can be given in terms of Thurston laminations as defined in \cite{ThurstonLaminations}). Then the two rays at angles $\theta/2$ and $(1+\theta)/2$ land at the critical point, so the critical point is biaccessible and the two angles $\theta/2$ and $(1+\theta)/2$ are biaccessible angles (with respect to the angle $\theta$). We say that $\pair{\theta/2}{(1+\theta)/2}$ forms the \emph{critical angle pair} (compare Figure~\ref{FigDefKneading}). Similarly, all further angles on the backwards orbit of $\theta$ are biaccessible angles. By induction on $k$ (the number of iterations required to reach $\theta$), each angle $\phi$ on the backwards orbit of $\theta$ has a unique angle $\phi'\neq \phi$ on the backwards orbit of $\theta$ with $\step(\phi)=\step(\phi')$ that is not separated by precritical angle pairs with lower values of $\step$, and then $\pair{\phi}{\phi'}$ forms a precritical angle pair. In particular, all precritical points are biaccessible, which forms a countable set. Therefore, to find the Hausdorff dimension of biaccessible angles in the dynamics modeled by the external angle $\theta$, we only need to investigate angles $\phi$ that are not on the backwards orbit of $\theta$, or equivalently we only need to investigate itineraries in $\Symg=\{\0,\1\}^\Nplus$. 

In parameter space, it is known that all parameter rays at periodic angles of fixed period $n>1$ land in pairs \cite{Orsay,MiOrbits,ExtRayMandel}, so their landing points are biaccessible and all periodic angles are biaccessible angles in parameter space. If two parameter rays at periodic angles $\theta_1$ and $\theta_2$ (necessarily of equal period) land together, we say that $\pair{\theta_1}{\theta_2}$ forms a \emph{periodic parameter angle pair}. Since all periodic external angles are thus biaccessible (in parameter space, at least when the period is $2$ or greater), and these are exactly the angles with $\*$-periodic kneading sequences and they form a countable set, it follows that the Hausdorff dimension of biaccessible angles in parameter space is determined by angles with kneading sequences in $\Symg$, and periodic parameter angle pairs give a necessary condition, and this condition turns out to be sufficient except for sets of dimension zero).

We can now give a combinatorial definition of biaccessibility, both in the dynamical plane and in parameter space.

\begin{defin}{(Combinatorial Biaccessibility and Angle Pairs).}
\label{Def:CombBiac} \lineclear 
An angle $\varphi \in \Circle$ is called 
\emph{combinatorially biaccessible with respect to $\theta\in\Circle$}
if there is a $\varphi' \in \Circle$ with $\varphi'\neq\varphi$ such that no precritical angle pair separates $\phi$ from $\phi'$. 
We call $\pair{\varphi}{\varphi'}$ a {\em dynamic angle pair}. 

An angle  $\theta \in \Circle$ with non-periodic kneading sequence in $\Symg$ is called  \emph{combinatorially biaccessible} (in parameter space)
if there is a $\theta' \in \Circle$ with $\theta'\neq\theta$ so that no periodic parameter angle pair $\pair{\theta_1}{\theta_2}$ separates $\theta$ from $\theta'$.  We call $\pair{\theta}{\theta'}$ a \emph{parameter angle pair}.
\end{defin}

\begin{rem}
All precritical angles are combinatorially 
biaccessible according to this definition.
Indeed, if $\theta$ is not combinatorially biaccessible, 
then $\pair{\theta/2}{(\theta+1)/2}$
is an dynamic angle pair that is not separated by any other precritical angle pair,
and this carries over to preimage angle pairs of $\pair{\theta/2}{(\theta+1)/2}$.
If, on the other hand, $\theta$ is combinatorially biaccessible, 
say $\pair{\theta}{\theta'}$ forms the corresponding dynamic angle pair, then 
 $\pair{\theta/2}{\theta'/2}$ and  $\pair{(\theta+1)/2}{(\theta'+1)/2}$
are angle pairs that are not separated by precritical angle pairs.
It is thus not relevant that
 $\pair{\theta/2}{(\theta+1)/2}$ and $\pair{\theta'/2}{(\theta'+1)/2}$
separate each other (while landing at the same point). 
\end{rem}

Our combinatorial estimates will all be with respect to this definition (in dynamical planes and in parameter space). In order to explain how this definition relates to the topological concept of biaccessibility as introduced earlier, we start with a simple lemma.

\begin{lem}{(Combinatorially Biaccessible Angle Pairs).} 
\label{Lem:CombBiac}\lineclear
An angle $\phi$ is combinatorially biaccessible with respect to $\theta$ if and only if there is a $\phi'\in\Circle$ with $\phi'\neq\phi$ so that $\nu_\theta(\phi)=\nu_\theta(\phi')$.
\end{lem}

\begin{proof}
If two angles $\phi\neq\phi'$ are separated by precritical angle pairs, then the angle pair with lowest value of $\step$ is always unique (because two angle pairs with equal value of $\step$ are always separated by an angle pair with lower value, which is easily confirmed inductively).

If two angles $\phi\neq\phi'$ are separated by a unique precritical angle pair $\pair{\theta_1}{\theta_2}$ with $\step(\theta_1)=\step(\theta_2)=k$, then $\nu_\theta(\phi)$ and $\nu_\theta(\phi')$ differ in the $k$-th position. 

Conversely, if the $k$-th entries of $\nu_\theta(\phi)$ and $\nu_\theta(\phi')$ 
are different, then $D^{\circ(k-1)}(\phi)$ and $D^{\circ(k-1)}(\phi')$ are separated by the diameter $\pair{\theta/2}{(1+\theta)/2}$, and by taking $k-1$ preimages it follows that $\phi$ and $\phi'$ are separated by a precritical angle pair $\pair{\theta_1}{\theta_2}$ with $\step(\theta_1)=\step(\theta_2)=k$.
\end{proof}

The corresponding statement in parameter space will be discussed below.  

\subsection{Combinatorial and topological biaccessibility}

We will now relate combinatorial biaccessibility with respect to $\theta$ to topological biaccessibility in a Julia set for which the dynamic ray at angle $\theta$ lands at the critical value. 

\begin{prop}{(Topologically and Combinatorially Biaccessible Angles).}
\label{Prop:ExtraAngles}\lineclear
For every quadratic polynomial with connected Julia set, the set of topologically biaccessible angles is a subset of the set of combinatorially biaccessible angles. In every case, these two sets have the same Hausdorff dimension.
\end{prop} 
\begin{proof}
If two rays at angles $\phi$ and $\phi'$ land together at a point that is not on the backwards orbit of the critical point, then $\phi$ and $\phi'$ clearly cannot be separated by a precritical angle pair. 
Therefore, combinatorial biaccessibility in the sense of Definition~\ref{Def:CombBiac} is necessary for topological biaccessibility in a dynamical plane. 

The converse is true at least in locally connected Julia sets; however, some care is necessary because local connectivity in itself does not seem to imply easily that combinatorially biaccessible angle pairs actually have rays that land at a common point. However, the somewhat stronger property of trivial fibers is sufficient, see \cite{Fibers1,Fibers2}: by definition, a point $z$ within a Julia set $\J$ has trivial fiber if, for every $z'\in \J\sm\{z\}$, there is a pair of dynamic rays at periodic angles that separates $z'$ from $z$. This implies local connectivity of $\J$ at $z$, but it is locally a strictly stronger property. In particular, it implies that every $z'\in \J\sm\{z\}$ can also be separated from $z$ by a pair of rays that land at the same precritical point. Hence, if all fibers are trivial, then every angle $\phi$ that is combinatorially biaccessible is in fact topologically biaccessible: there is another angle $\phi'$ so that the rays at angles $\phi$ and $\phi'$ land together. (Note that in the case of bounded Fatou components, a slight modification is necessary in the definition of fibers, allowing for separation lines to pass through bounded Fatou components \cite{Fibers1}.)

In fact, all currently known proofs of local connectivity of Julia sets of quadratic polynomials show that all fibers are trivial: this is true for polynomials with attracting or parabolic periodic orbits (Douady and Hubbard~\cite{Orsay}), as well as for all quadratic polynomials for which all periodic orbits are repelling and that are not infinitely renormalizable (Yoccoz \cite{HubbardYoccoz}; for the renormalizable case see also \cite{Fibers3}), and for various infinitely renormalizable quadratics in work of Kahn and Lyubich (see \cite{LyubichQuadratics,KahnLyubich,KahnLyubich2} and the references therein), and finally in the case of a Siegel disk of bounded type \cite{Petersen}. A separate argument is the general result that if a polynomial Julia set is locally connected at every point, then all fibers are trivial \cite{Fibers3} (as mentioned, this implication is not true at individual points $z$).

The relation between combinatorial and topological biaccessibility may thus be non-trivial only in the presence of irrationally indifferent cycles or for infinitely renormalizable polynomials. In fact, for a quadratic polynomial with an indifferent cycle of period $1$, the set of angles that is not separated from the indifferent fixed point by a precritical angle pair has Hausdorff dimension zero \cite{BulSen} (this is the set of external angles of the Siegel disk in case the Julia set is locally connected). All other angles cannot be part of a biaccessible angle pair \cite{SchleicherZakeri}, so biaccessibility concerns only a set of Hausdorff dimension zero (including the backwards orbit of all rays that are not separated from the indifferent fixed point, which is a countable union), so this does not affect our dimension estimates.

If there is an indifferent periodic cycle of period $M>1$, then the set of external angles related to this cycle still has dimension zero, and all other points have trivial fibers \cite{HubbardYoccoz,Fibers3}, so combinatorially and topologically biaccessible angles coincide.

Finally, for infinitely renormalizable polynomials
(see Definition~\ref{DefRenormalizable}), we note first that, for the same reason as before, all points $z\in \J$ that are outside of the little $M$-renormalizable Julia set for some $M$ have trivial fibers, while external angles corresponding to $M$-renormalization are in a set of Hausdorff dimension at most $1/M$ (see Section~\ref{SecRenorm}). Therefore, if a polynomial is infinitely renormalizable, then all angles belong to trivial fibers with the exception of a set of dimension zero, which again does not affect our estimates.
\end{proof}

\begin{prop}{(Topologically and Combinatorially Biaccessible Angles in $\M$).}
\label{Prop:TopoCombBiaccM}\lineclear
For the Mandelbrot set, the set of topologically biaccessible angles is a subset of the set of combinatorially biaccessible angles, and both sets have the same Hausdorff dimension.
\end{prop}

\begin{proof}
If two parameter rays land together, then they cannot be separated by a parameter ray-pair at periodic angles (the landing point of the latter is not the landing point of any parameter ray except the two periodic ones). This implies that an angle can only be topologically biaccessible if it is combinatorially biaccessible in the sense of Definition~\ref{Def:CombBiac}. 

We will now argue that ``most'' combinatorially biaccessible angles are topologically biaccessible (except possibly for a set of dimension zero). Of course, the set of external angles that do not land has always Hausdorff dimension zero, but this is not a sufficient argument: we also need to consider those rays that land, but possibly not at the same point as another ray at a combinatorially associated angle.

The set of angles in parameter space that are associated to $M$-renormalizable parameters has again Hausdorff dimension at most $1/M$ (again Section~\ref{SecRenorm}), so the set of angles that are associated to infinitely renormalizable parameters has dimension zero and does not affect our estimates. 
Similarly, the set of angles that are associated to any particular hyperbolic component have Hausdorff dimension zero, see Corollary~\ref{CorHausdorffParaAngles}. 

Therefore, except for a set of angles of Hausdorff dimension zero, every parameter ray is associated to a fiber of $\M$ corresponding to polynomials for which all cycles are repelling and which are not infinitely renormalizable, and such fibers are trivial by \cite{HubbardYoccoz}. Fibers of $\M$ are defined in terms of parameter ray-pairs at periodic angles, so if two parameter rays are not separated by parameter ray-pairs of periodic angles (and are not in the exceptional set), then they belong at the same fiber and thus land at the same point; so for all angles except for a set of dimension zero, combinatorial biaccessibility implies topological biaccessibility. 
\end{proof}

\begin{rem}
In parameter space, it may be worth noting that parameter rays at periodic angles (which have $\*$-periodic kneading sequences) always land in pairs and are thus topologically biaccessible. Moreover, there are uncountably many further parameter rays that accumulate at hyperbolic components: their external angles have periodic (but not $\*$-periodic) kneading sequences, and these rays land at a boundary point of a hyperbolic component (with irrationally indifferent dynamics). These rays are all not separated by periodic parameter ray-pairs and thus combinatorially biaccessible, but they are the only rays landing at the same point, so they are not topologically biaccessible. These angles (from a zero-dimensional set,
see Corollary~\ref{CorHausdorffParaAngles}) thus form a difference between the two definitions of biaccessibility. A more precise definition of combinatorial biaccessibility would thus be ``an angle $\phi$ is combinatorially biaccessible if either there is an angle $\phi'$ so that $\phi$ and $\phi'$ are not separated by a periodic angle pair, or if the kneading sequence $\nu(\phi)$ is periodic but not $\*$-periodic''. This definition would be more cumbersome without strengthening Proposition~\ref{Prop:TopoCombBiaccM} --- however, if all fibers of $\M$ were trivial (which is equivalent to local connectivity of $\M$ \cite{Fibers2}), then topological and combinatorial biaccessibility of $\M$ in this sense would coincide.
\end{rem}

Definition~\ref{Def:CombBiac} is thus a good combinatorial description of biaccessibility, ignoring topological subtleties.

\subsection{Combinatorial biaccessibility, itineraries, and kneading sequences}

In order to investigate biaccessibility from a combinatorial point of view, we translate it to the setting of itineraries. In order to do this, we extend the definition of the $\rho$-function to itineraries.

For a kneading sequence $\nu \in \Syms$ and $x=\nu_\theta(\varphi)\in \Symg$, let
\begin{equation*}\label{eq:rhox}
\rhox(n) := \min\{ k > n \colon x_k \neq \nu_{k-n}\}.
\end{equation*}
Obviously, $\rho_{\nu, \nu}=\rho_\nu$ for $\rho_\nu$ as in
Definition~\ref{DefRho}.

\begin{lem}{(Condition for Itinerary to be Biaccessible).}
\label{LemBiaccess} \lineclear
Let $\theta,\varphi\in\Circle$ be two arbitrary external angles, let $\nu:=\nu(\theta)$ be the kneading sequence of $\theta$ and let $x:=\nu_\theta(\phi)$ be the itinerary of $\varphi$ with respect to $\theta$; assume that $x\in\Symg$.  

Then $\varphi$ is combinatorially biaccessible with respect to $\theta$ if and only if there is a $k\ge 2$ that satisfies
\begin{equation*}\label{EqBiaccess}
\orb_{\rhox}(1) \cap \orb_{\rhox}(k) = \emptyset.
\end{equation*}
\end{lem}

\begin{proof}
We use symbolic dynamics modeled after the situation that the critical value is the landing point of the ray at angle $\theta$. Then $\pair{\theta/2}{(1+\theta)/2}$ is the critical angle pair with $\step=1$, and all other precritical angle pairs are preimages of this one. We say that a precritical angle pair $\pair{\theta_k}{\theta'_k}$ with $\step(\theta_k)=\step(\theta'_k)$ is a \emph{closest precritical angle pair} (to $\phi$) if it is not separated from $\phi$ by any precritical angle pair for which the value of $\step$ is less than $\step(\theta_k)$. Since any two precritical angle pairs with equal value of $\step$ are separated by another angle pair with lower value, for any value of $\step$ there can be at most one closest precritical angle pair.

The closest precritical angle pair $\pair{\theta_1}{\theta'_1}$ with 
$\step(\theta_1)=\step(\theta'_n)=1$ is always $\pair{\theta/2}{(1+\theta)/2}$. 
Now define a sequence $(m_j)$ with $m_1=1$ and, if $\pair{\theta_j}{\theta'_j}$ is a closest precritical angle pair with $\step(\theta_j)=\step(\theta'_j)=m_j$, then let $\pair{\theta_{j+1}}{\theta'_{j+1}}$ be the precritical angle pair with lowest value of $\step$ that separates $\phi$ from $\pair{\theta_j}{\theta'_j}$;  we define $m_{j+1}:=\step(\theta_{j+1})=\step(\theta'_{j+1})$. 

This $\rho$-itinerary 
describes the itinerary of $\theta_j$ and of $\theta'_j$ (which is the same) completely: 
the first $m_j-1$ entries coincide with that of $\phi$, \ie with $\nu_\theta(\phi)$, 
then comes a $\star$, and then the itinerary of $\theta$ (which is $\nu(\theta)$). 
This means $\rho_{\nu,x}(m_j)=m_{j+1}$ for all $j$, so $(m_j)=\orb_{\rho_\nu,x}(1)$ (equality of sequences).

Now suppose $\pair{\phi}{\phi'}$ is a biaccessible angle pair; by definition, this means that no closest precritical angle pair $\pair{\theta_j}{\theta'_j}$ separates $\phi$ from $\phi'$, so they all (for $j>1$) separate both $\phi$ and $\phi'$ from the angle pair $\pair{\theta/2}{(1+\theta)/2}$. 

Without loss of generality, we may assume that $0<\phi<\phi'<1$. Then there is a precritical angle $\psi_1\in(\phi,\phi')$ and hence a precritical angle pair $\pair{\psi_1}{\psi'_1}\subset(\phi,\phi')$; set $k_1:=\step(\psi_1)=\step(\psi_2)$.

As before, for a precritical angle pair $\pair{\psi_j}{\psi'_j}$,
define a sequence $k_j$ so that $\pair{\psi_{j+1}}{\psi'_{j+1}}$ is the precritical angle pair with lowest value of $\step$ that separates $\phi$ from $\pair{\psi_j}{\psi_{j+1}}$; we define $k_{j+1}:=\step(\psi_{j+1})=\step(\psi'_{j+1})$. Then all $\{{\psi_j},{\psi'_j}\}\subset(\phi,\phi')$, so none of the $\pair{\psi_j}{\psi'_j}$ separate $\pair{\phi}{\phi'}$ from $\pair{\theta_1}{\theta'_1}$. Therefore, the sequences $(m_j)$ and $(k_j)$ are disjoint. 

Finally, comparing itineraries, it is easy to check that $(k_j)=\orb_{\rho_{\nu,x}}(k_1)$, and this proves the existence of two disjoint $\rho_{\nu,x}$-orbits as required. 

The converse is similar: $\orb_{\rho_{\nu,x}}(1)$ always describes the sequence of closest precritical angle pairs $\pair{\theta_j}{\theta'_j}$, starting at the diameter, that separate the previous angle pair from $\phi$. If we label so that $\theta_j<\theta'_j$, then the sequence $(\theta_j)$ is monotonically increasing and $(\theta'_j)$ is monotonically decreasing. Denote their limits by $\theta_\infty\le\theta'_\infty$. Observe that if $\theta_\infty<\theta'_\infty$, then $\pair{\theta_\infty}{\theta'_\infty}$ forms a biaccessible angle pair. Even if $\phi\not\in\{\theta_\infty,\theta'_\infty\}$, both $\pair{\phi}{\theta_\infty}$ and $\pair{\phi}{\theta'_\infty}$ are biaccessible angle pairs. In order to prove that $\phi$ is combinatorially biaccessible, it thus suffices to prove that $\theta_\infty<\theta'_\infty$.

This is assured by the existence of  $\orb_{\rho_{\nu,x}}(k)$: this describes a similar sequence starting with the closest precritical angle pair for which $\step$ has the value $k$. If these sequences are disjoint for some $k$, then none of the associated angle pairs can be separated from $\phi$ by any $\pair{\theta_j}{\theta'_j}$, and every angle pair described by $\orb_{\rho_{\nu,x}}(k)$ provides a uniform lower bound for $\theta'_j-\theta_j$, so their limits are different as required.
\end{proof}

Now we turn to parameter space.

\begin{lem}{(Biaccessible Kneading Sequence).}
\label{LemNuBiaccess} \lineclear
An angle $\theta$ is combinatorially biaccessible (in parameter space) in the sense of Definition~\ref{Def:CombBiac} if and only if its  kneading sequence
$\nu = \nu(\theta)$ has a $k\ge 2$ such that
\begin{equation*}\label{EqBiaccessPara}
\orb_{\rho}(1) \cap \orb_{\rho}(k) = \emptyset.
\end{equation*}
\end{lem}

\begin{figure}[htp]
\includegraphics{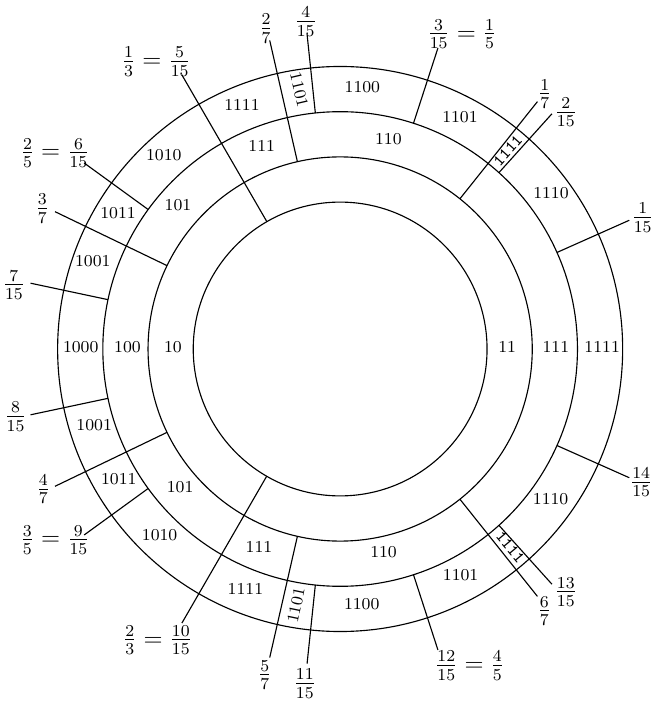}
\caption{Illustration of the map $\theta \mapsto \nu(\theta)$.
The angles of the form $a/(2^k-1)$, $k \le n$, divide $\Circle$ into arcs 
that map to components of $n$-cylinder sets.}\label{FigNu}
\end{figure}

\begin{proof}
The structure of this proof is similar to Lemma~\ref{LemBiaccess}; we use the structure given by periodic parameter angle pairs $\pair{\psi}{\psi'}$ of equal period that land together at a common point. 

Note first that, as a function of the angle $\psi$, the $n$-th entry in the kneading sequence $\nu(\psi)$ changes between $\0$ and $\1$ exactly when $\psi$ is periodic of period $n$ or dividing $n$. Since angles of equal period come in pairs, it follows that if two non-periodic angles $\theta$ and $\theta'$ are not separated by a periodic angle pair, their kneading sequences must coincide; and if the periodic ray pair of lowest period that separates them has period $n$, then the two kneading sequences differ for the first time at the $n$-th position  (this statement uses the fact, known as Lavaurs' Lemma \cite{Lavaurs}, 
that if any two parameter angle pairs $\pair{\psi_1}{\psi'_1}$ and $\pair{\psi_2}{\psi'_2}$ of equal period are such that $\pair{\psi_1}{\psi'_1}$ separates $\pair{\psi_2}{\psi'_2}$ from the angle $0$, then the two angle pairs are separated from each other by an angle pair of lower period).

Conversely, if their kneading sequences (which must be in $\Symg$) differ for the first time at the $n$-th position, then the lowest period of an angle pair separating $\theta$ and $\theta'$ is exactly $n$. In particular, if the two kneading sequences coincide, then they are not separated by any periodic ray pair.

Now suppose that a non-periodic angle $\theta$ is combinatorially biaccessible in parameter space, \ie there is an angle $\theta'$ that is not separated from $\theta$ by any periodic parameter angle pair. Then $\nu(\theta)=\nu(\theta')=:\nu$. The angle pair $\pair{\theta_1}{\theta'_1}$ of lowest period $n_1$ that separates the angle $0$ from $\pair{\theta}{\theta'}$  must have period $n_1=\rho_\nu(1)$, and inductively the angle pair of lowest period that separates $\pair{\theta_j}{\theta'_j}$ from $\pair{\theta}{\theta'}$, say $\pair{\theta_{j+1}}{\theta'_{j+1}}$, has period $n_{j+1}=\rho_\nu(n_j)$. 
The sequence of periods of these angle pairs is thus $\orb_{\rho_\nu}(1)$; this is the internal address of $\theta$. 

Sorting the angle pairs again such that $0<\theta<\theta'<1$ and $0<\theta_j<\theta'_j<1$ for all $j$, it follows that the sequence $(\theta_j)$ is monotonically increasing and converges to a limit in $(0,\theta]$, while $(\theta'_j)$ is monotonically decreasing and converges to a limit in $[\theta',1)$. 

Now pick a periodic angle $\psi_1\in(\theta,\theta')$ of least possible period and denote this period by $k_1$. Then the angle pair of lowest period separating $\psi_1$ from $\theta$, and thus also from $\theta'$, has period $k_2:=\rho_\nu(k_1)$, and by induction we obtain a sequence of angle pairs $\pair{\psi_j}{\psi'_j}$ so that $\pair{\psi_{j+1}}{\psi'_{j+1}}$ separates $\pair{\psi_j}{\psi'_j}$ from $\theta$ and from $\theta'$, and their periods are $\orb_{\rho_\nu}(k_1)$. All their angles are in $(\theta,\theta')$, so these angle pairs are disjoint from the $\pair{\theta_j}{\theta'_j}$. Again by Lavaurs' Lemma, all periods $n_j$ and $k_{j'}$ must be different. Therefore, $\orb_{\rho_\nu}(1) \cap \orb_{\rho_\nu}(k_1) = \emptyset$.

\looseness-1
For the converse, consider an  angle $\theta$ with kneading sequence $\nu=\nu(\theta)$ and a $k\in\Nplus$ such that $\orb_{\rho_\nu}(1) \cap \orb_{\rho_\nu}(k) = \emptyset$. If $\theta$ is periodic, then it is part of a ray pair and thus biaccessible, so we may assume that $\theta$ is not periodic.
Construct the sequence of periodic angle pairs $\pair{\theta_j}{\theta'_j}$ as above so that each of these ray pairs separates its successor from the angle $0$. Let $\theta_\infty$ and $\theta'_\infty$ be the two limits of their angles; they satisfy $\theta_\infty\le\theta'_\infty$. We want to show that $\theta_\infty<\theta'_\infty$. By construction, none of the angles $\theta$, $\theta_\infty$, and $\theta'_\infty$ can be separated from each other by a periodic parameter angle pair, and while possibly $\theta\in\{\theta_\infty,\theta'_\infty\}$, it follows in any case that $\theta$ is combinatorially biaccessible if $\theta_\infty<\theta'_\infty$. 

To complete this proof, we need to find an angle $\psi$ with $\theta_\infty<\psi<\theta'_\infty$. More precisely, we will find a $k'\in\orb_{\rho_\nu}(k)$ so that the kneading sequence $\nu^{(k')}$ consisting of the first $k'-1$ entries in $\nu$, followed by $\*$ and continued periodically, is admissible, and so that there is an angle $\psi_{k'}$ with $\nu(\psi)=\nu^{(k')}$; this angle $\psi_{k'}$ can then be chosen in $(\theta_\infty,\theta'_\infty)$. Not every such $\nu^{(k')}$ is admissible; see the remark after the proof. However, we have the following.

\begin{claim}{}\label{Claim:Admissible}
If $k' \in \orb_{\rho_\nu}(k)$ is such that $k' > \rho_\nu(m)$
for all $m < \rho_\nu(k)$, then there exists a periodic angle $\psi_{k'}\in\Circle$ of period $k'$ with $\nu(\psi_{k'})=\ovl{\nu_1 \dots \nu_{k'-1}\star}$.
\end{claim}

We first complete the proof of the lemma using this claim. 
There may be several choices for $\psi_{k'}$; choose one for which the angled internal address coincides longest possible (among all candidates) with $\theta$. It will turn out that this determines $\psi_{k'}$ uniquely, but if there are still several choices, then the choice is arbitrary. 

The angle $\psi_{k'}$ must be part of a ray pair $\pair{\psi_{k'}}{\psi'_{k'}}$. Inductively let $\pair{\psi_{j+1}}{\psi'_{j+1}}$ be the ray pair of least period that separates $\pair{\psi_{j}}{\psi'_{j}}$ from $\theta$ (such a ray pair always exists because $\theta$ is non-periodic, so the kneading sequences never coincide). We do not know yet that $\psi_{k'}$ or $\psi_j$ are in $(\theta_\infty,\theta'_\infty)$.

We claim that the ray pair $\pair{\psi_{j+1}}{\psi'_{j+1}}$ is unique. If not, then two ray pairs of equal period separate $\pair{\psi_{j}}{\psi'_{j}}$ from $\theta$, say $\pair{\psi_{j+1}}{\psi'_{j+1}}$ and $\pair{\tilde\psi_{j+1}}{\tilde \psi'_{j+1}}$. Again by Lavaurs' lemma, neither can separate the other from the angle $0$, and this is possible only if (possibly by switching labels) $\pair{\psi_{j+1}}{\psi'_{j+1}}$ separates $\pair{\psi_{j}}{\psi'_{j}}$ from $0$, while $\pair{\tilde\psi_{j+1}}{\tilde\psi'_{j+1}}$ separates $\theta$ from $0$, and both are not separated from each other by a ray pair of lower period, so they have equal kneading sequences and thus equal internal addresses. By \cite[Theorem~4.3]{IntAddr} there exists another ray pair $\pair{\tilde\psi_{k'}}{\tilde\psi'_{k'}}$ of equal period and with equal kneading sequence as $\pair{\psi_{k'}}{\psi'_{k'}}$, but separated from $0$ by $\pair{\tilde\psi_{j+1}}{\tilde\psi'_{j+1}}$. This contradicts the choice of $\psi_{k'}$ as the one for which the angled internal address coincides for the longest possible time with that of $\theta$. Therefore, all ray pairs $\pair{\psi_{j+1}}{\psi'_{j+1}}$ are uniquely determined by $\psi_{k'}$. 

Our next claim is that no $\pair{\psi_j}{\psi'_j}$ separates $\theta$ from the origin. Indeed, if it does, then the sequences $\pair{\psi_j}{\psi'_j}$ and $\pair{\theta_{j'}}{\theta'_{j'}}$ will eventually coincide, but their periods are $\orb_{\rho(\nu)}(k)$ and $\orb_{\rho(\nu)}(1)$ and these are disjoint by hypothesis. 

Therefore, all $\pair{\psi_j}{\psi'_j}$ separate $\psi_{k'}$ and hence $\pair{\psi_{j-1}}{\psi'_{j-1}}$ from $\theta$ and from $0$; this implies that $\theta\not\in (\psi_j,\psi'_j)\supset(\psi_{j-1},\psi'_{j-1})$ for all $j$. Moreover, we must have $(\psi_j,\psi'_j)\subset(\theta_{j'},\theta'_{j'})$ for all $j$ and $j'$: otherwise, there are $j$ and $j'$ for which $(\psi_j,\psi'_j)\cap (\theta_{j'},\theta'_{j'})=\emptyset$, and then eventually $\pair{\theta_{j'}}{\theta'_{j'}}$ must appear among the $\pair{\psi_j}{\psi'_j}$, and this is an impossibility. 

Together, we have $\psi_{k'} \in (\psi_j,\psi'_j)\subset \bigcap_{j'}(\theta_{j'},\theta'_{j'})=(\theta_\infty,\theta'_\infty)$, which proves the lemma.
\end{proof}

We still need to prove Claim~\ref{Claim:Admissible}; for easier reference, we repeat the statement in a self-contained way.
\emph{If $\nu\in\Sigma^1$ is an admissible kneading sequence, $k\in\N$ satisfies $\orb_\rho(k)\cap\orb_\nu(1)=\emptyset$, and $k' \in \orb_\nu(k)$ is such that $k' > \rho(m)$
for all $m < \rho(k)$, then there exists a periodic angle $\psi_{k'}\in\Circle$ of period $k'$ with $\nu(\psi_{k'})=\ovl{\nu_1 \dots \nu_{k'-1}\star}$.}
\medskip

\begin{proof}[Proof of the claim]
Recall that $\nu$ (with internal address $1 \IntAddr S_1 \IntAddr S_2 \IntAddr \dots$) is admissible. The kneading sequence $\nu':=\ovl{\nu_1 \dots \nu_{k'-1}\star}$ is $\*$-periodic of period $k'$ and thus has an associated Hubbard tree $(T,f)$, where $f$ is a continuous self-map of the finite tree $T$ in which the unique critical point $c_0$ has period $k'$ and kneading sequence $\nu'$ \cite[Theorem~2.5]{BKS}.

The main task is to prove that $\nu'$ does not fail the admissibility condition in Definition~\ref{DefAdmissCond}. This implies, by \cite[Theorem~4.2]{BS1}, that the tree $T$ can be embedded into $\C$ so that the dynamics $f$ can be extended continuously to a neighborhood of $T$ and in fact to all of $\C$, and thus it is the Hubbard tree of a complex quadratic polynomial. The critical value has two characteristic external rays, and their external angles then have period $k'$ and kneading sequence $\nu'$. Every tree $(T,f)$ embedded into $\C$ gives thus rise to two angles as required in the claim, and there may be different choices for the embedding of the tree.

In order to give a proof by contradiction, suppose that $\nu'$ fails the admissibility condition (Definition~\ref{DefAdmissCond}) for some period $m_*$. This depends only on the first $\rho(m_*)$ entries of $\nu'$. Since $\nu$ is admissible and coincides with $\nu'$ for $k'-1$ entries, the non-admissibility of $\nu'$ cannot be determined by looking only at its first $k'-1$ entries. This implies $\rho_\nu(m_*) \geq k'$.
Since the endpoints of $T$ lie on the critical orbit, and there are more
endpoints than branch points in a tree, we have $m_*< k'$ 
and therefore $\rho_{\nu_\star}(m_*)=k'$. 
By hypothesis of the claim, $\rho_\nu(m)<k'$ for all $m < \rho_\nu(k)$, and hence $m_* \ge \rho_\nu(k)$.

The fact that $\nu'$ fails the admissibility condition for period $m_*$ means that the Hubbard tree $(T,f)$ has an $m_*$-periodic evil periodic orbit 
with $q \ge 3$ arms at each of its points \cite[Lemma 3.6]{BS1}. Denote the characteristic point of this evil orbit be $z_1$. Then the global arms at $z_1$ can be labeled $G_0,G_1,\dots,G_{q-1}$ so that 
$G_0 \owns c_0$ (the critical point) and so that after $m_*$ iterations, $G_1$ maps homeomorphically onto $G_2$, $G_2$ maps homeomorphically onto $G_3$, etc., $G_{q-2}$ maps homeomorphically onto $G_{q-1}$, and $G_{q-1}$ maps to its image that intersects both $G_0$ and $G_1$ (while $G_0$ contains the critical point, so its image contains $G_1$ and much (if not all) of $T$). This implies that the itineraries of $z_1$ and $c_1$ coincide for at least $(q-2)m_*$ entries.

Within the Hubbard tree, take $k-1$ consecutive preimages of the critical point $c_0$, always choosing the branch so that the itinerary of the resulting point, to be called $\zeta_k$, with respect to $c_0$ starts with $\nu_1\nu_2\dots\nu_{k-1}\star$. Since the critical value has period $k'>k$, there is no ambiguity. It may of course happen that the appropriate preimage is not in $T$; in this case, it is straightforward to extend $(T,f)$ appropriately (as in \cite[Lemma~3.6]{BS1}). The extended tree, say $T'$, comes with a continuous self map $f'$ that extends $f$ on $T$, and it shares all axioms of $(T,f)$ except minimality; in particular, $T'\sm\{c_0\}$ has at most two connected components and $f'$ is injective on each. 

We claim that $c_1$ separates $\zeta_k$ from $0$ and $z_1$, so that ``$\zeta_k$ is behind $c_1$'' (as seen from $c_0$). To show this, first observe that $k'\not\in\orb_\nu(1)$ (because $k' \in \orb_\nu(k)$ and $\orb_\nu(k)\cap\orb_\nu(1)=\emptyset$). In this case, the cutting time argument applied to $[c_0,c_1]$ shows that if we replace $c_1$ by a pair of points, one with entry $\0$ and one with $\1$ instead of the $\*$, then the point on the side of $c_0$ is such that it does not generate an entry $k'$, while $\zeta_k$ does generate such an entry ($k'\in\orb_\nu(k)$), so $\zeta_k$ is on the opposite side of $c_1$ than $c_0$, as claimed.

But this implies that $\zeta_k$ is in the same global arm of $z_1$ as $c_1$ and thus survives at least $(q-2)m_*\ge m_*$ iterations homeomorphically without hitting the critical point, which is in contradiction to $k<\rho_\nu(k)\le m_*$. 
\end{proof}

\begin{rem}
It is not true that if a kneading sequence $\nu$ is admissible, then for every $k$ the $\*$-periodic sequence of period $k$ that coincides with $\nu$ for $k-1$ entries is also admissible, even when $\orb_\rho(k)\cap\orb_\rho(1)=\emptyset$.
A counter-example is supplied by the kneading sequence corresponding to internal address 
$1 \IntAddr 2\IntAddr 4 \IntAddr 6 \IntAddr 7 \IntAddr 9 \IntAddr 11
 \IntAddr 13 \IntAddr 15 \IntAddr 17 \dots$ and $k=10$.
Then $\orb_{\rho}(k) \cap \orb_\rho(1) = \emptyset$, but
$1 \IntAddr 2\IntAddr 4 \IntAddr 6 \IntAddr 7 \IntAddr 9 \IntAddr 10$
is not admissible (evil period $m_* = 5$).
\end{rem}

\section{Preliminaries on Cylinder Sets} \label{sec:CoreEntropy}

We write $C_{e_1e_2\dots e_n}=\{x_1x_2\dots\in\Symg\colon 
x_1=e_1,\dots, x_n=e_n\}$ for cylinder sets in $\Symg$ of length $n$. 
Denote the length of a cylinder $C$ by $|C|$, that is: $|C_{e_1e_2\dots e_n}|=n$.

\begin{lem}{(The Shape of the $\nu_\theta$-Inverse of Cylinders).}
\label{LemNuThetaPreimage} \lineclear
The map $\nu_\theta: \Circle \to \{ \0 , \star, \1\}^{\Nplus}$
is (in general) non-injective: for each $n$-cylinder
$C_{e_1\dots e_n}\in\{\0,\1\}^n$, the preimage
$\nu_\theta^{-1}(C_{e_1\dots e_n})$ consists of at most $n$ open arcs of total length
$2^{-n}$. 
\end{lem}

\begin{rem}
There is no a priori lower bound (in terms of $n$) 
on the length of the components
of $\nu_\theta^{-1}(C_{e_1\dots e_n})$.
Indeed, such components have endpoints $\varphi$ and $\varphi'$ 
(not necessarily in that order) satisfying
$2^k\varphi - a = \theta = 2^m\varphi' - b$ for integers
$a, b \in \N$ and $k \le m \le n$.
If $k = m$, then $a \neq b$ and $|\varphi'-\varphi| \ge 2^{-m}$, 
but otherwise
\[
|\varphi'-\varphi| \ge 2^{-m} d( (2^{m-k}-1) \theta, \Z),
\] 
and for Liouville numbers $\theta$ this lower bound can be extremely small 
compared to $2^{-n}$. 
\end{rem}

\begin{proof}
The open arcs $A_\0 = (\frac{\theta}{2}, \frac{\theta+1}{2})$ and
$A_\1 = (\frac{\theta+1}{2}, \frac{\theta}{2})$ form the partition
of $\Circle$ which yields itineraries in $\{ \0, \1 \}^\Nplus$.
(Thus we ignore those countably many $\varphi \in \Circle$ whose 
itinerary $\nu_\theta(\varphi)$ contains a $\star$.)
 
The two sets $A_\0$ and $A_\1$ correspond to the two $1$-cylinders $C_\0$ 
and $C_\1$ of $\Symg$. 
Suppose by induction on $n$ that the set $A$ corresponding to the $n$-cylinder
$C_{e_1\dots e_n}$ has at most $n$ components $I$.
For $j \in \{ \0, \1 \}$, $D^{-1}(I) \cap A_j$ consists of one interval, or
two if $\theta \in I$.
Therefore $D^{-1}(I) \cap A_j$ which corresponds to the $n+1$-cylinder
$C_{e_1\dots e_n j}$ has at most $n+1$ components.
This proves the induction step and hence the lemma.
\end{proof}

\begin{lem}{(The Shape of the $\nu$-Inverse of Cylinders).}
\label{LemNuPreimage} \lineclear
The map $\nu: \Circle \to \Syms$
is non-injective: for each $n$-cylinder
$C_{e_1\dots e_n}$, the preimage
$\nu^{-1}(C_{e_1\dots e_n})$ consists of finitely many open arcs 
of length between $2^{-(n+1)}$ and $2^{-(2n+1)}$.
\end{lem}

\begin{rem}  
The total number of arc-components of $\nu(C_{e_1\dots e_n})$ is based on
an estimate in how many ways we can embed a periodic Hubbard tree into 
the plane, see Lemma~\ref{LemEmbedding}. We estimate this number
in Lemma~\ref{LemNumberEmbed2}.
\end{rem}

\begin{proof}
This time the arc components of $\nu^{-1}(C_{e_1\dots e_n})$ are open arcs
with endpoints $\theta$ and $\theta'$ satisfying
$2^{k+1} \theta - a = \theta$, $2^{m+1} \theta' - b = \theta'$, that is
$\theta = a/(2^{k+1}-1)$, $\theta' = b/(2^{m+1}-1)$ for some
$k \le m \le n$ and $a,b \in \N$, see Figure~\ref{FigNu}.
Taking $k=m=n$ and $|a-b| = 1$ we get the upper bound
$|\theta-\theta'| \le 1/(2^{n+1}-1)$.
The lower bound is
$\min \{  |a/(2^{k+1}-1) -  b/(2^{m+1}-1)| > 0 \ : \ a,b \in \N\} 
\ge 1/(2^{2n+1}-1)$.
\end{proof}

\section{Preliminaries on Hausdorff Dimension} \label{sec:Hausdorff}

The motor for the dimension estimates will be the following
elementary lemma.

\begin{lem}{(Hausdorff Dimension of Sample Sets).}
\label{LemBasic_hdim}\lineclear
Given integers $u>v \geq 1$,
construct nested compact sets $A_s\subset[0,1]$ (for $s\geq 0$) as follows:
\begin{itemize}
\item
Let $A_0 = [0,1]$;
\item
Divide each of the $(u-v)^s$ intervals of $A_s$ of length
$u^{-s}$ into $u$ equal intervals and remove the closures of $v$
of them, chosen arbitrarily. Then $A_{s+1}$ is the union of the
closures of the remaining $(u-v)^{s+1}$ intervals of length
$u^{-(s+1)}$ each.
\end{itemize}
Let $A = \cap_s A_s$.
Then
\(\displaystyle \hdim(A) = \frac{\log(u-v)}{\log u} \).
\end{lem}

\begin{proof}
Since $A_s$ consists of $(u-v)^s$ intervals of length $u^{-s}$, the
box dimension of $A$ is $\frac{\log(u-v)}{\log u}$.
Therefore $\hdim(A) \leq \frac{\log( u-v)}{\log u}$.

For the lower bound we use the measure $\mu_s$ on $A_s$ which
assigns mass $(u-v)^{-s}$ to each of the $(u-v)^s$ intervals of
$A_s$ and refine it to a measure $\mu$ on $A$ using Kolmogorov's
extension theorem, see \eg \cite{Chung}. For a point $x\in A$, let $I_s(x)$
be the interval of
$A_s$ containing $x$. If $B(x;\eps)$ denotes the $\eps$-ball around $x$,
then
the interval $I_s(x)$ is contained in the ball $B(x;\eps)$
for $u^{-s} < \eps \leq u^{-s+1}$, and $B(x;\eps)$ is contained in at most
$u$ intervals of length $u^{-s}$.
Choosing $\delta < \frac{\log u-v }{\log u }$, we have
for $s$ sufficiently large:
$$
\mu(B(x;\eps)) \leq u(u-v)^{-s} = u \left(u^{-s \cdot \frac{\log u-v 
}{ \log u } }\right)
\leq \eps^{-\delta}.
$$
The Frostman Lemma (see \eg \cite{Mattila})
now implies that the $\delta$-dimensional Hausdorff mass of $A$ is positive.
Since  $\delta < \frac{\log u-v }{\log u }$ is arbitrary, we obtain the
required lower bound $\hdim(A) \geq \frac{\log u-v }{ \log u }$.
\end{proof}

As usual, we endow $\Symg = \{ \0 , \1 \}^{\Nplus}$ with the metric
$d(x, \tilde x)=\sum_{i\geq 1}|x_i-\tilde x_i|2^{-i}$.
The \emph{binary representation map}
$\binary \colon \Symg \to \Circle=\R/\Z$ is given
by $\binary(x_1x_2\ldots)=\sum x_i 2^{-i}$; it is injective
except for the countably many dyadic rationals. We define the Hausdorff dimension of $Y\subset \Symg$ by $\hdim(Y):=\hdim(\binary(Y))$, and we denote the Hausdorff dimension of subsets of $\Circle$  and of $\Symg$ by $\hdim$.

\begin{cor}{(Hausdorff Dimension of Concatenations of Blocks).}
\label{CorBasic_hdim}\lineclear
For distinct blocks $X_1, \dots , X_k$ of $\0$s and $\1$s, none of 
which is a suffix of another, let
$$
B = \left\{ x = W_1W_2\dots\ :\  W_i \in \{ X_1, \dots , X_k\}
\text{ for all } i \makehigh \right\} \subset \Symg 
$$
(in other words, consider an arbitrary infinite concatenation of blocks $X_i$).  
Then
$\displaystyle \hdim(B) \geq \frac{\log k}{m \log 2}$ for $m = \max_i |X_i|$.
\end{cor}

\begin{proof}
Extend each block $X_i$ to the left to a block $\tilde X_i$ of length 
$m$ in an arbitrary way. Since no $X_i$ is a suffix
of any other $X_j$, the resulting blocks $\tilde X_i$ are distinct.
Then Lemma~\ref{LemBasic_hdim} immediately gives
that $\hdim(\tilde B) = \log k/(m \log 2)$ for $\tilde B = \{ x = 
W_1W_2\dots\ :\ W_i
\in \{ \tilde X_1, \dots , \tilde X_k\} \}$.
Indeed, $\tilde B$ can be transformed into a subset of $\Circle$
using the binary extension map $\binary:\Symg \to \Circle$ that makes 
the shift on $\Symg$ commute
with the angle doubling map on $\Circle$.

Now define $h:\Symg \to \Symg$ by replacing every `non-overlapping'
occurrence
of a block $X_i$ in $x = x_1x_2x_3\dots \in \Symg$ by the block 
$\tilde X_i$ and leaving
the other coordinates $x_j$ untouched. More precisely, we work from 
left to right:
whenever we encounter a block $X_j$ not overlapping with an 
occurrence of some block $X_i$
replaced previously, then we replace it with $\tilde X_j$.
Then $h$ maps $B$ bijectively and continuously onto $\tilde B$, and 
the Lipschitz constant
of $h$ is at most $1$. Therefore $\hdim(B) \geq \hdim(\tilde B)  = 
\log k/(m \log 2)$, as
required.
\end{proof}

\begin{lem}{(Symbolic Codings that Preserve Hausdorff Dimension).}
\label{LemHdimPreserved}\lineclear
Let $P:{\Nplus} \to {\Nplus}$ be a polynomial and $K > 0$ and suppose that
$I: \Circle \to  \Symg$ is a map such that the preimage
$I^{-1}(C)$ of any $n$-cylinder consists
of at most $P(n)$ intervals of length $\leq K2^{-n}$.
Then $\hdim(I^{-1}(\Omega)) \leq \hdim(\Omega)$  for any set $\Omega 
\subset  \Symg$.

If $Q:{\Nplus} \to {\Nplus}$ is a polynomial such that
the preimage $I^{-1}(C)$ of any $n$-cylinder contains an arc 
of length $\geq 2^{-n}/Q(n)$,
then $\hdim(I^{-1}(\Omega)) \geq \hdim(\Omega)$  for any set $\Omega 
\subset \Symg$.
\end{lem}

\begin{proof}
Let $\eps > 0$ be arbitrary and take any $\delta'' > \delta' > \delta
= \hdim(\Omega)$.
Let $N$ be so large that
\begin{itemize}
\item $K^{\delta''} P(n) <
2^{n(\delta''-\delta')}$ for all $n \geq N$;
\item  $\sum_i \diam(U_i)^{\delta'} < \eps$, where $\{ U_i \}$ is a
cover of $\Omega$ such that $\diam(U_i) < 2^{-N}$ for each $i$.
(Without loss of generality we can assume
that each $U_i$ is a cylinder set of length $n_i\geq N$.
In the standard metric on $\Sym$, $\diam(U) = 2^{-|U|}$.)
\end{itemize}
Then $\{I^{-1}(U_i) \}_i$ defines a countable cover $\{ V_j\}_j$ of
$I^{-1}(\Omega)$, each interval $V_j$ has length at most $K 2^{-n_i}$, and
\begin{eqnarray*}
\sum_j \diam(V_j)^{\delta''}
&=& \sum_{n \geq N} \sum_{|U_i| = n}
\sum_{V_j \subset I^{-1}(U_i)} \diam(V_j)^{\delta''} \leq
\sum_{n \geq N} \sum_{|U_i| = n}  P(n) K^{\delta''}2^{-n\delta''} \\
&\leq&
\sum_{n \geq N} \sum_{|U_i| = n} P(n)K^{\delta''}
2^{-n(\delta''-\delta')} \diam(U_i)^{\delta'}\\
&\leq&
\sum_{n \geq N} \sum_{|U_i| = n} \diam(U_i)^{\delta'} < \eps\; .
\end{eqnarray*}
Since this is true for every $\eps > 0$ and $\delta'' > \delta$,
it follows that $\hdim(I^{-1}(\Omega)) \leq \delta$.

Now for the second statement, take $0 < \delta'' < \delta' < \delta
= \dim_H(\Omega)$ and $K > 0$ arbitrary.
Then there exists $N$ so large that
\begin{itemize}
\item $( \frac1{Q(n)} 2^{-n} )^{\delta''} > 2^{-\delta' n}$ for all $n
\geq N$;
\item $\sum_i (2^{ -|C_i| } )^{\delta'} > 2K$,
where $\{ C_i \}$ is a cover of $\Omega$
with cylinder sets with $|C_i| \geq N$.
\end{itemize}
For each $C_i$, let $A_i$ be an interval in $I^{-1}(C_i)$
of length $\geq 2^{-n}/Q(n)$.
Let $\{ V_j \}_j$ be any open cover of $I^{-1}(\Omega)$
with intervals of length $< 2^{-n}/(2Q(n))$.
For each $i$, let ${\mathcal V}_i = \{ V_j : V_j \subset A_i\}$,
so $\sum_{V_j \in {\mathcal V}_i} |V_j|^{\delta''} > \frac12
|A_i|^{\delta''}$.
Therefore
\begin{eqnarray*}
\sum_j |V_j|^{\delta''} &\geq&
\sum_i \sum_{V_j \in {\mathcal V}_i} |V_j|^{\delta''}
\geq \frac12 \sum_i |A_i|^{\delta''} \\
&\geq& \frac12 \sum_i \left(\frac{1}{Q(|C_i|)} 2^{-|C_i|} \right)^{\delta''}
\geq \frac12 \sum_i 2^{-\delta' |C_i|} > K.
\end{eqnarray*}
Since $K$ and $\delta'' < \delta' < \delta$ are arbitrary, we obtain
$\hdim(I_{\theta}^{-1}(\Omega)) \geq \delta$.
This proves the lemma.
\end{proof}

\section{Dimension for combinatorially biaccessible 
itineraries}\label{sec:comb_dim_sequences}

We will produce two pairs of bounds for the Hausdorff dimension of
biaccessible itineraries (and kneading sequences). These constitute the main step for proving Theorems~\ref{ThmHdimBiaccAngles} and~\ref{ThmHdim_mandel2}.

\begin{prop}{(Dimension of Biaccessible Sequences).}
\label{PropHdimBiacSymbolics}\lineclear
(i) For any kneading sequence $\nu$, the Hausdorff dimension of
biaccessible itineraries with respect to $\nu$ is in $\I(N(\nu),
\kappa(\nu))$.

\noindent
(ii) The Hausdorff dimension of biaccessible kneading
sequences $\nu$ with $N(\nu) = N$ and $\kappa(\nu) = \kappa$
is in $\I(N, \kappa)$.
\end{prop}

\begin{proof}
(i)
Fix a kneading sequence $\nu$, let $N=N(\nu)$ and define, for $k\geq 2$,
\[
    \Biack = \left\{ x \in  \Symg \colon  k =
\min\{ i\geq 2 \colon \orb_{\rhox}(1) \cap \orb_{\rhox}(i) = \emptyset\}
\makehigh \right\}.
\]
Note that the set of biaccessible itineraries is $\bigcup_k\Biack$. 
We will show that all $\Biack$ satisfy the same dimension bounds.

\heading{Upper Bound \boldmath $U_1(N)$:\unboldmath} In order to prove that
$\hdim(\Biack) \leq \frac{\log(2^N-1)}{\log 2^N}$, we show that for 
sufficiently large $n>k$, every $n$-cylinder $C_{e_1\dots e_n}$ 
contains at least one $n+N$-cylinder that is disjoint from $\Biack$.

Choose $x \in C_{e_1\dots e_n}$.
Let $a := \max \{ i \leq n\colon i \in \orb_{\rhox}(1) \}$
and $b := \max \{ i \leq n\colon i \in \orb_{\rhox}(k) \}$.
Then clearly $\rhox(a)>n$ and $\rhox(b)>n$.
Suppose that $C_{e_1\dots e_n}\cap \Biack\neq\emptyset$ (otherwise 
there is nothing to show); then
$a\neq b$.
Let $w_a = \nu_{n-a+1}\dots \nu_{n-a+N}$ and $w_b =
\nu_{n-b+1}\dots \nu_{n-b+N}$. Recall that $\rhox(a)$ finds the
first difference between
$x_{a+1}x_{a+2}\dots$ and $\nu_1\nu_2\dots$. Therefore
$\rhox(a)\leq n+N$ unless $x$ starts with $e_1\dots 
e_a\nu_1\nu_2\dots\nu_{n-a} w_a$,
and similarly for $b$.

Our task is the following: given $\nu$, $n$ and $e_1\dots e_n$, we
want to find at least one $n+N$-cylinder in $C_{e_1\dots e_n}$
  disjoint from $\Biack$.

Let $\0\dots \0$ and $\0\1\dots \0$ be the two words of length $N$
which contain no $\1$, except possibly at the second position.
The following three cases are easy to check:

\begin{description}
\item
[Case 1] $w_a \neq \0\dots \0$ and $w_b \neq \0\dots \0$.
We claim that the cylinder $C_{e_1\dots e_n\0\dots \0}$ is disjoint
from $\Biack$. Indeed, for $x\in C_{e_1\dots e_n\0\dots \0}$, we
get $\rhox(a)\in\{n+1,\dots,n+N\}$: we have $\rho_{\nu,x}(a)\geq n+1$ 
by definition of $a$, and $\rho_{\nu,x}(a)\leq n+N$ because 
$w_a\neq\0\dots\0$ means that $\nu$ does not have a sequence of $N$ 
zeroes starting at position $n-a+1$. After 
$\rho_{\nu,x}(a)\in\{n+1,n+N\}$, the orbit $\orb_{\rhox}(a)$ 
increases in steps of $1$ until it reaches
$n+N$, hence $n+N\in\orb_{\rhox}(a)\subset\orb_{\rhox}(1)$.
Similarly, $n+N\in\orb_{\rhox}(k)$, which proves the claim.

\item
[Case 2] $w_a = \0\dots \0$ and $w_b = \1 \dots $.
This time, we claim that $C_{e_1\dots e_n\0\1\dots \0}$ is
disjoint from $\Biack$: for $x\in C_{e_1\dots e_n\0\1\dots \0}$,
we have $\rhox(a)=n+2$, and after that, $\orb_{\rhox}(a)$
increases in steps of $1$ up to $n+N$, so again
$n+N\in\orb_{\rhox}(a)$. This time, $\rhox(b)=n+1$ and
$\rhox(\rhox(b))=n+N$, so
$n+N\in\orb_{\rhox}(1)\cap\orb_{\rhox}(k)$.

\item
[Case 3] $w_a = \0\dots \0$ and $w_b = \0 \dots $.
Now the entire $N+1$-cylinder $C_{e_1\dots e_n\1}$ is disjoint
from $\Biack$: for $x\in C_{e_1\dots e_n\1}$, we have
$\rho(a)=n+1=\rho(b)$.
\end{description}
These three cases cover all possibilities, possibly interchanging the 
roles of $w_a$ and $w_b$.
Hence each $n$-cylinder contains at least one $n+N$-cylinder
disjoint from $\Biack$.
By Lemma~\ref{LemBasic_hdim}, $\Biack$ is contained in a Cantor set of
Hausdorff dimension 
\[
\frac{\log(2^N-1)}{\log 2^N} = 1+ \frac{\log(1-2^{-N})}{\log 2^N} \le 1-\frac{2^{-N}}{\log 2^N}
\]
as claimed (using the standard bound $\log(1+x)\le x$ for $|x|<1$).

\heading{Lower Bound \boldmath $L_1(N)$:\unboldmath}
First assume that $N = N(\nu) \geq 6$.
Observe that the beginning of $\nu = \1\0\0\dots\0\1\dots\,\,$
contains $N-3$ zeroes in a row. Take $M = \lfloor N/2 \rfloor -1$,
and let
\begin{equation}\label{eq:BW}
B = \left\{ x =  W_1W_2\dots \in \Symg \ : \ 
W_i\in\{0,1\}^M\sm\{\0\0\dots\0\}  \right\}.
\end{equation}
In particular, no $x \in B$ contains $N-3$ consecutive symbols $\0$.
It follows that $\rhox(i)< i+ N$ for all $i \geq 1$. Moreover,  every $x
\in B$ admits 
at least two disjoint $\rhox$-orbits: if ${m_i}$ describes the 
positions of the entries $\1$ in $x$, then 
$\rho_{\nu,x}(m_i-1)=m_{i+1}$, and then the orbit increases in steps 
of $1$ until $m_{i+2}-1$, and it later reaches $m_{i+4-1}-1$ etc. A 
different orbit goes through $m_{i+1}-1$, $m_{i+3}-1$, $m_{i+5-i}$, 
etc., and is disjoint from the first one. Therefore, all sequences in 
$B$ are biaccessible.

According to Lemma~\ref{LemBasic_hdim}, the Hausdorff dimension of $B$ is
\begin{align*}
\frac{\log(2^{\lfloor N/2\rfloor -1}-1)}{\log 2^{\lfloor N/2 \rfloor -1}}
&= 1+\frac{\log(1-2^{-\lfloor N/2-1\rfloor })}{\log 2^{\lfloor N/2 \rfloor -1}}
\ge 1-\frac{2^{-\lfloor N/2-1\rfloor}}{(1-2^{-\lfloor N/2-1\rfloor})\log 2^{\lfloor N/2 \rfloor -1}}
\\
&=1-\frac{1}{(2^{\lfloor N/2-1\rfloor}-1)\log 2^{\lfloor N/2 \rfloor -1}}
\end{align*}
(using the standard lower estimate $\log(1+x)\ge x/(1+x)$ for $|x|<1$). 

Now let us treat the case $N = 5$, so $\nu = \1\0\0\1\dots$.
In this case, we take
\[
B = \{ x =  W_1W_2\dots \in \Symg \ : \  W_i \in\{ \1\1,\1\0\} \text{ 
for } i \geq 2 \}.
\]
Then every $x \in B$ admits two disjoint $\rhox$-orbits.
The Hausdorff dimension of $B$ is
$\log 2 / \log 4 = 1/2$ according to Lemma~\ref{LemBasic_hdim}.
This proves the lower bound $L_1(N)$.

\begin{rem}
The same idea gives lower bounds for other beginnings of kneading sequences:
\begin{align*}
\nu &= \1\0\1\1\0\ldots:\ \text{ taking } W_i = \1\1\1\1 \text{ or } 
\1\0\1\0
\text{ gives } \hdim(B) \ge 1/4. \\
\nu &= \1\0\1\0\0\ldots:\ \text{ taking } W_i = \1\1\1\1\1 \text{ or 
}  \1\1\0\1\0
\text{ gives } \hdim(B) \ge 1/5. \\
\nu &= \1\0\1\1\1\1\ldots: \text{ taking } W_i = \1\0\1\1\1\0 \text{ 
or }  \1\1\1\0\1\0
\text{ gives } \hdim(B) \ge 1/6.
\end{align*}
Incidentally, for the latter two examples, these bounds equal the
respective bounds $L_2(S_\kappa)$ below.
The bound $1/4$ for $\nu =  \1\0\1\1\0\dots$ is better than 
$L_2(S_\kappa) = 1/5$.
\end{rem}

\heading{Upper Bound \boldmath $U_2(S_\kappa)$: \unboldmath}
We start with Case (c) in the definition of $L_2(S_\kappa)$ and $U_2(S_\kappa)$; see \eqref{Eq_U2kappa}; in this case $S_\kappa < S_{\kappa+1} < \infty$.
By the definition of $\kappa$, we can write $S_j = p_{j-1}S_{j-1}$ for
$1 \leq j \leq \kappa$ and we define $p_\kappa := \max\{ i \geq 1 : iS_\kappa < S_{\kappa+1}\}$.
Then $S_{\kappa+1} \leq 2p_\kappa S_\kappa$ and
\begin{equation}\label{Eq_form_nu}
\nu_1\dots \nu_{S_j} =
(\nu_1\dots \nu_{S_{j-1}})^{p_j-1}(\nu_1\dots \nu'_{S_{j-1}}),
\quad \nu_1\dots \nu_{p_\kappa S_\kappa} = (\nu_1\dots \nu_{S_\kappa})^{p_\kappa},
\end{equation}
where $\nu'_i = \1$ if $\nu_i = \0$ and vice versa.
Every $n < p_\kappa S_\kappa$ can be written uniquely as
$$
n = \sum_{j=0}^{\kappa} a_j S_j, \qquad 0 \leq a_j < p_j.
$$
If $\rho(n) < S_{\kappa+1}$, then
\begin{equation}\label{Eq_decomp}
\rho(n) = n + aS_h \text{ for } h \geq \min\{ j : a_j \neq 0\}
\text{ and some } 1 \leq  a < p_h.
\end{equation}
Now if $x \in \Biack$, then we can enumerate the entries
of $\orb_{\rho_{\nu,x}}(1)$ and $\orb_{\rho_{\nu,x}}(k)$
  as $1 = u_0 < u_1 < \dots$ and  $k = v_0 < v_1 < \dots$.
We try to estimate how many different
sequences $(u_s)_{s \geq 0}$ and  $(v_t)_{t \geq 0}$ (and hence 
sequences $x \in \Biack$) can be both disjoint and satisfy
$\rho_{\nu,x}(u_s) = u_{s+1}$ and $\rho_{\nu,x}(v_t) = v_{t+1}$.

If $x$ is known up to entry $u_{s+1}$ and $u_s < v_t < u_{s+1}$,
then $v_{t+1}$ is fully determined, provided $v_{t+1} < u_{s+1}$.
Let us analyze what can happen if
  $u_s < v_t < u_{s+1} < v_{t+1}$.

{\bf Claim:} If $u_{s+1}-u_s \leq p_\kappa S_\kappa$, then
$u_{s+1}-v_t = aS_h$ for some $h < \kappa$ and $1 \leq  a < p_h$,
and furthermore $v_{t+1}-v_t \geq \min\{ S_{h+1}, S_\kappa + 1\}$.

To prove this, let $y \in \Sym$ be such that the first difference 
between $x$ and $y$ is at position $u_{s+1}$.
Abbreviate $n = v_t-u_s$ and $r = u_{s+1}-u_s$.
Then $y_{u_s+1} \dots y_{u_{s+1}} = \nu_1 \dots \nu_r$,
and $\rho_{\nu,y}(v_t) = u_{s+1}$. By \eqref{Eq_decomp}, this means that
$u_{s+1}-v_t = r-n = aS_h$ for some $h$ and $1 \leq a < p_h$.
Now the form of $\nu$ given in \eqref{Eq_form_nu} shows that
$v_{t+1} - v_t = \rho_{\nu,x}(v_t) - v_t \geq S_{h+1}$, because
$\rho_{\nu,x}(v_t) - v_t < S_{h+1}$ would imply
that $\rho_{\nu,x}(v_t) = \rho_{\nu,x}(u_{s+1})$.
This proves the claim.

Next take $s'$ maximal such that $u_{s'} < v_{t+1}$.
Then one of the following holds:
\begin{enumerate}
\item $v_{t+1}-v_t = a'S_j$ for some $h+1 \leq j \leq \kappa$ and
$1 \leq a' < p_j$. Then \eqref{Eq_form_nu} shows that $v_{t+1}-u_{s'} =
a'' S_{j-1}$ for some $1 \leq a'' < p_{j-1}$, and the above argument
(with the roles of $u$ and $v$ interchanged) implies that
$u_{s'+1}-u_{s'} \geq S_j$.
\item  $v_{t+1}-v_t > S_{\kappa}$. In fact, if $S_{\kappa+1} > 
v_{t+1}-v_t > S_{\kappa}$ then using  \eqref{Eq_form_nu} again, we 
find $v_{t+1} = u_{s'+1}$, so
in this case  $v_{t+1}-v_t \geq S_{\kappa+1}$, but of course
$S_{\kappa+1} = S_\kappa + 1$ is possible.
\end{enumerate}
Let us say that $u$ and $v$ {\em switch roles at entries $(s,t)$} if
$u_s < v_t < u_{s+1} < v_{t+1}$. Let $h = h(s,t)$ be such that 
$u_{s+1}-v_t = aS_h$.
The above arguments show that if $h(s,t) \leq \kappa$ and
$v$ and $u$ switch roles again
at entries $(t+1, s')$, then $h(t+1, s') \geq h(s,t)$.
Hence, at switches, $h$ is non-decreasing at least until it
exceeds $\kappa$, whereas
between switches (say $u_s < v_t < v_{t+1} < u_{s+1}$, the entry
$v_{t+1}$ is fully determined by $u_{s+1}$.
To illustrate this, let us give an example:
$$
\nu = \1 \0 \1\0\1\1 \1\0\1\0\1\0 \1\0\0\dots
$$
with $\kappa = 3$ and internal address
$1 \IntAddr 2\IntAddr 6 \IntAddr 12 \IntAddr 15 \IntAddr \dots$,
and
\[
x = \underbrace{\1}_{ u_0=1}
\overbrace{\1}^{ v_0,\ h = 0}
\overbrace{\0}^{ v_1,\ h = 1} \1\ \ \0\ \ \1
\overbrace{\1}^{v_2,\ h = 2}\1\ \ \0\ \ \1\ \ \0\ \ \1
  \underbrace{\1}_{ u_1, \ h=2}\1\ \ \0\ \ \1\ \ \0\ \ \1
\overbrace{\1}^{\tiny v_3,\ h = 2}\dots
\]
with $k = 2$.
We see that $h$ stays constant if the roles of $u$ and $v$ switch,
and increases between switches.
Furthermore, two consecutive switches of roles takes $S_{h+1}$ digits.

We can code the consecutive switches by integers
$l_j \geq 0$: each $l_j$ indicates the number of pairs of switches
where $h(s,t) = h(s+1,t+1)$
(or $h(t,s) = h(t+1,s+1)$ for reversed roles of $u$ and $v$)
remain constant at $j$. If $l_j = 0$, it means that $h(s,t)$
increases from below $j$ to above $j$. Suppose there are $r \geq 1$ occurrences of
$h(s,t) > \kappa$ before $h(s,t)$ drops to $\leq \kappa$ again.
Let $m_j \geq S_\kappa+1$, $1 \leq j \leq r$,
denote the distances between
the remaining switches before $h(s,t) \leq \kappa$ again.
(If $r=0$, then there are no such $m_j$s.)
Thus the whole loop from $h(s,t)=0$ to the last $h(s,t) \geq \kappa$
takes at least
$\sum_{j=0}^\kappa l_j S_{j+1} + \sum_{j = 1}^r m_j$ digits.
Let us introduce a second index $i$ to indicate the loop number.
Then the pair $(l_{i,j})_{i=1, j=0}^{n,\ \kappa}, (m_{i,j})_{i=1, 
j=1}^{n,\ r_i}$
encodes
a cylinder set in $\Biack$ going through $n$ loops, and the cylinder length
is at least
$k + \sum_{i=1, j = 0}^{n, \kappa} l_{i,j} S_{j+1} + \sum_{i=1, j = 
1}^{n, r_i} m_{i,j}$.

Let $\delta = U_2(S_\kappa)$.
The cylinders encoded by  $(l_{i,j})_{i=1, j=0}^{n,\ \kappa}, 
(m_{i,j})_{i=1, j=1}^{n,\ r_i}$
form a cover of $\Biack$ with diameter $< 2^{-(k+nS_\kappa)}$.
Its $\delta$-dimensional Hausdorff measure is bounded by
$$
2^{-k\delta} \sum 2^{-\delta\ \sum_{i=1, j=0}^{n, \kappa} l_{i,j} 
S_{j+1}} \cdot
\sum 2^{-\delta\ \sum_{i=1, r_i = 1, j=1}^{n, \infty, r_i} m_{i,j} },
$$
where the first main sum runs over all combinations of $n(\kappa+1)$
positive integers $l_{i,j}$ and the second main sum over all
combinations of integers $m_{i,j} \geq p_\kappa S_\kappa+1$.
Using geometric series, the estimate $\sum_{l=0}^\infty 2^{-l\alpha}
\leq 1+\int_0^\infty 2^{-x\alpha} dx = 1+\frac{1}{\alpha \log 2}$,
and changing the order of product and sum, we can rewrite this quantity as
\begin{align} \label{eq:dimHquantity}
& 2^{-k\delta} \cdot \prod_{i=1}^n \left[
\prod_{j=0}^{\kappa} \sum_{l_{i,j} = 0}^\infty 2^{-\delta l_{i,j} 
S_{j+1}} \cdot
\sum_{r_i=1}^\infty \prod_{j=1}^{r_i} \sum_{m_{i,j} = p_\kappa S_\kappa+1}^\infty
2^{-\delta m_{i,j}} \makehigh \right] \nonumber \\
& \qquad \leq \  2^{-k\delta} \cdot
\prod_{i=1}^n \left[
\prod_{j=0}^{\kappa} \left( 1+ \frac{1}{\delta S_{j+1} \log 2} \right)
\cdot \sum_{r_i=1}^\infty \left( \frac{2^{-\delta (p_\kappa S_\kappa+1)} 
}{1-2^{-\delta}}
\right)^{r_i} \makehigh \right] .
\end{align}
Next observe that $S_{\kappa+1} \geq p_\kappa S_\kappa+1 > S_j \geq 2^j$ for $0 \leq j  \leq \kappa$, and hence
$\prod_{j=0}^\kappa (1+ \frac{1}{\delta S_{j+1} \log 2})
\leq 2^{1/\delta (\log 2)^2}$.
The second factor is another geometric series, and can be computed as
$$
\sum_{r_i=1}^\infty \left( \frac{2^{-\delta (p_\kappa S_\kappa+1)} 
}{1-2^{-\delta}} \right)^{r_i}
= 2^{-\delta (p_\kappa S_\kappa+1) - \log(1-2^{-\delta} - 2^{-\delta 
(p_\kappa S_\kappa+1)})/\log 2}.
$$
Therefore expression \eqref{eq:dimHquantity}
is bounded by $2$ to the power
$$
-k\delta + \frac{n}{\log 2}
\left( \frac{1}{\delta \log 2} - 
 (p_\kappa S_\kappa+1)\ \delta \log 2
- \log(1-2^{-\delta} - 
2^{-\delta (p_\kappa S_\kappa+1)}) \right),
$$
so that \eqref{eq:dimHquantity} is bounded in $n$ if and only if
the factor in the brackets above is non-positive.
In coordinates $\delta \log 2 = a/\sqrt{p_\kappa S_\kappa+1}$, this
is equivalent to
\begin{equation}\label{eq:P}
P(a, p_\kappa S_\kappa) := e^{a\sqrt{p_\kappa S_\kappa+1}} 
(1-e^{-a/\sqrt{p_\kappa S_\kappa+1}} ) - e^{\sqrt{p_\kappa S_\kappa+1}/a} \geq 1.
\end{equation}
For $p_\kappa S_\kappa = 3$, this can be solved numerically by
$a \geq a_0 \approx 1.8234 \log 2 > 1.2638$. Since $P(a, p_\kappa S_\kappa)$
is increasing in $a$ and $p_\kappa S_\kappa$ for $a \geq a_0$,
it follows that \eqref{eq:dimHquantity} is bounded in $n$
for $\delta = \frac{1.8234}{\sqrt{p_\kappa S_\kappa + 1}}$.
Since $p_\kappa S_\kappa+1 > S_{\kappa+1}/2$ and
$1.8234 \sqrt{2} < \sqrt{7}$, we get
$\hdim(\Biack) \leq U_2(S_\kappa)$.

\begin{rem}
If $S_\kappa = 3$, then $\nu = \1\1\0\dots$ and
$\delta = \frac{1.8234}{\sqrt{p_\kappa S_\kappa + 1}} = 0.9117$ gives a slightly
better estimate than
$U_1(N) = U_1(3) = \log 7/\log 8 \approx 0.9358$.
If $S_\kappa = 2$, then $\nu = \1\0\dots$ and $U_1(N) \leq U_1(4)$
is the better upper bound.
At the other end, given any $\eps > 0$, we can take
$\delta = \frac{1+\eps}{\sqrt{p_\kappa S_\kappa+1} \log 2}$
as upper bound provided $p_\kappa S_\kappa$ is sufficiently large.
\end{rem}

Cases (a) and (b) are limit cases, and enforce the sequences
$v_{t+1}-v_t$ and $u_{s+1}-u_s$ to be non-decreasing for every $x \in E_k$.
Therefore for any $N$ and co-countably many $x \in E_k$, eventually
$v_{t+1}-v_t \geq N$ and $u_{s+1}-u_s \geq N$. This means that eventually,
there are at most two ``free choices'' of symbol in $x$ within every
$N$ entries. Hence we can find $c>0$ such the number of $n$-cylinders 
in $E_k$ is at most $c2^{2n/N}$ for every $n \in \N$,
whence $\hdim(\Biack) \leq 2/N$. But $N$ is arbitrary, so the upper bound $U_2(S_\kappa) = 0$ holds in these cases too.

\heading{Lower Bound \boldmath $L_2(S_\kappa)$:\unboldmath}
Write $\nu = \nu_1\nu_2\dots$ and define
$V = \nu_1\nu_2\dots  \nu_{S_{\kappa}-1} \nu'_{S_{\kappa}}$ and
$\hat V = \nu_1\nu_2\dots \nu_{S_{\kappa+1}-1} \nu'_{S_{\kappa+1}}$,
where $\nu'_i = \1$ if $\nu_i = \0$ and vice versa.
(Note that if $S_{\kappa+1} = \infty$, then there is nothing to prove.)
Let
\begin{equation}\label{eq:BVV}
B = \left\{ x =  W_1W_2\dots \in \Symg \ : \  W_i \in\{ V,\hat V\} \right\}.
\end{equation}
Corollary~\ref{CorBasic_hdim} implies
$
\hdim(B) \geq \frac{1}{S_{\kappa+1}}
$
as claimed, so it suffices to show that each $x \in B$ admits
two disjoint $\rhox$-orbits.
By construction of $x \in B$,
$\rhox^{\circ i}(S_\kappa) = |W_1W_2\dots W_i|$ for all $i \geq 0$.
We will show that the $\rhox$-orbit of $S_{\kappa-1}$ is disjoint from
this.
Note that $V$ is the concatenation of $S_\kappa / S_{\kappa-1}$ blocks
$\nu_1\nu_2\dots \nu_{S_{\kappa-1}}$.
Therefore, for any integer $a\in\{1,\dots,
\frac{S_\kappa}{S_{\kappa-1}}-1\}$,
\[
\rho_{\nu, VV}(aS_{\kappa-1}) = \rho_{\nu, V\hat V}(aS_{\kappa-1})
= S_\kappa + aS_{\kappa-1},
\]
where we extended the definition of $\rhox$ to the case where $x$ is a
finite block in the obvious way.
Also
\[
\rho_{\nu, \hat VV}(aS_{\kappa-1}) = \rho_{\nu, \hat V\hat
V}(aS_{\kappa-1}) =
S_\kappa.
\]
Let $n = S_{\kappa+1} - S_\kappa$, so we can write
$\hat V = VW$ for $W = \nu_1\dots \nu_n$.
Furthermore $W = V^iX$ for some $i \geq 0$ and $m := |X| < S_\kappa$.
We can use \eqref{Eq_decomp} to compute $\rho_{XV, XV}(m) = aS_h$
for some $h \leq \kappa-1$ and $1 \leq a \leq S_{h+1}/S_h$.
If $h = \kappa-1$, then $X$ is the concatenation
of at most $p_{\kappa}-1$ blocks $\nu_1\dots \nu_{S_{\kappa-1}}$
and in this case we readily find
$$
\rho_{\nu, \hat VV}(S_{\kappa}) = \rho_{\nu, \hat V\hat V}(S_{\kappa}) =
S_\kappa + aS_{\kappa-1}.
$$
If $h \leq \kappa-2$, then
$$
\rho_{\nu, \hat VV}(S_{\kappa}) = \rho_{\nu, \hat V\hat V}(S_{\kappa}) =
S_\kappa + \rho(n) = S_{\kappa+1} + (\rho(n)-n) =  S_{\kappa+1} +
(\rho(m)-m).
$$
Since $\rho(m) \leq S_{\kappa-1}$, Lemma 4.2 in \cite{BKS} gives that
$S_{\kappa-1} \in \orb_\rho(\rho(m)-m)$, and therefore
$S_{\kappa+1} + S_{\kappa-1}$ belongs to the $\rho_{\nu, \hat VV}$-orbit
(and to the $\rho_{\nu, \hat VV}$-orbit) of $S_\kappa$.

Combining these facts, we derive that the $\rhox$-orbit of $S_{\kappa-1}$
contains  $|W_1\dots W_i| + a_iS_{\kappa-1}$ for each $i$
and some $1 \leq a_i < p_\kappa$,
and hence is disjoint from the $\rhox$-orbit of $S_\kappa$.
\medskip

(ii)
Now for the second statement, \ie for kneading sequences, we repeat 
the proof with
\begin{equation}\label{Eq_Symg_k}
\Symg_{N, S, k} = \left\{ \nu \in \Symg \ : \begin{array}{l}
N(\nu) = N,\ S_\kappa(\nu) = S,\ \\[1mm]
k = \min\{ i \ : \ \orb_{\rho}(1) \cap \orb_{\rho}(i) = \emptyset\}
\end{array}
\makehigh\right\}.
\end{equation}
Take $G = N$ or $S_{\kappa+1}$ according to whether
the dimension estimate is obtained from
\eqref{Eq_U1N} or \eqref{Eq_U2kappa}.
For $\nu \in \Symg_{N,\kappa, k}$, instead of comparing subwords of
$\nu$ with a fixed itinerary, we compare
subwords of $\nu$ with $\nu$ itself, and
in the above arguments, only a comparison with
$\nu_1 \dots \nu_G$ matters.
Therefore there is no change in the upper bounds, also if we have to exclude
the non-admissible kneading sequences.

For the lower bounds, take  $n > \max\{k , G\}$,
and we can always select an admissible $n$-cylinder
for $C$ intersecting $\Symg_{N,\kappa, k}$ (from equation~\ref{Eq_Symg_k}) and such that
$\orb_\rho(i) \owns n$ for all $i < G$.

If $G = N$, then the sequences
$B = \{ \nu = CW_1W_2\dots : W_i \in \{ \0, \1\}^M \sm \{ \0\0\dots\0\} \}$
(constructed in the same gist as \eqref{eq:BW})
have the property that $\rho(i)-i < G$ for all $i \geq n$, and hence
they satisfy Admissibility Condition~\ref{DefAdmissCond}.
If $G = S_{\kappa+1}$, then we use
$B = \{ \nu = CW_1W_2\dots : W_i \in \{ \0, \1\}^M \sm \{ V, \hat V \} \}$
(as in \eqref{eq:BVV}).
The same reasoning gives that all $\nu \in B$ satisfy the
  admissibility condition, and so we obtain the same lower bounds
$L_1(N)$ and $L_2(S)$.
\end{proof}

We can now translate the first (dynamical) half of 
Proposition~\ref{PropHdimBiacSymbolics} from itineraries to external 
angles of dynamic rays; the second (parameter) half with the transfer 
from kneading sequences to external angles of parameter rays will be 
treated in Section~\ref{Sec:DimParaAngles}.

\begin{proof}[Proof of Theorem~\ref{ThmHdimBiaccAngles}]
From Lemma~\ref{LemNuThetaPreimage} we know that for each $n$-cylinder
set $C_{e_1\dots e_n}$, the preimage
$\nu_\theta^{-1}(C_{e_1\dots e_n})$ consists of at most $n$ open arcs,
with combined length $2^{-n}$. Hence each of these arcs has length $\le 2^{-n}$
and at least one of them has length $\ge 2^{-n}/n$.

Using Lemma~\ref{LemHdimPreserved}, we can transfer the dimension bounds
of Propositions~\ref{PropHdimBiacSymbolics} 
to the combinatorially biaccessible dynamic angles,
proving the theorem.
\end{proof}

\begin{rem}
If $\nu$ is periodic (but not $\*$-periodic), this may correspond to a Siegel disk in the Julia set it models. There is  a Cantor set $K$ of dynamic angles with the same itinerary $\nu$. Lemma~\ref{LemNuThetaPreimage} doesn't fail:
it just says that $K$ can be covered by $n$ arcs of combined length $2^{-n}$
for each $n$ and hence $\dim_H(K) = 0$. This fact was already proved by Bullett and Sentenac \cite{BulSen}.
\end{rem}

\section{Dimension estimates for angles in parameter space}
\label{Sec:DimParaAngles}

In this section, we make the transition from the dimension of 
kneading sequences (Proposition~\ref{PropHdimBiacSymbolics} (ii)) to 
the dimension of external angles of rays in parameter space.

In \cite{BKS} we constructed Hubbard trees based on the combinatorial
information encoded in the internal address or kneading sequence only.
In \cite[Lemma 3.1.]{BS1} it was shown that all branch points that are not 
precritical have a
representative periodic point, called {\em characteristic point}
on their orbit that lies on the arc $[0, c_1]$ and closer to $c_1$
than any other periodic point on the same orbit.
The precise definition is as follows:

\begin{defin}{(Characteristic Point).}
\label{DefCharacteristic} \lineclear
A periodic point $p$ on a Hubbard tree is called {\em
characteristic} if $c_1$ lies in a different component of
$T \sm \{ p \}$ than every other of $\orb(p)$.
\end{defin}

Characteristic points come in two types, {\tame} and {\evil},
of which the {\tame} points are the ones that actually occur in
true embedded Julia sets. We call the components of $T \sm \{ z_1 \}$
the {\em global arms} of $z_1$, whereas the global arms intersected with
a small neighborhood of $z_1$ are called {\em local arms}.
The next lemma collects from  \cite[Lemma 3.6]{BS1} those properties
of global arms of branch points that are relevant for this paper.

\begin{lem}{(Global Arms at Branch Points Map Homeomorphically).}
\label{LemHomeo} \lineclear
Let $z_1$ be the characteristic point of a
{\tame} $n$-periodic orbit of branch points, each with $q \geq 3$ arms.
Then $n$ appears in the internal address, and
the global arms at $z_1$ can be labeled $G_0$, $G_1$, \ldots, $G_{q-1}$
so that $G_0 \owns 0$, $G_1 \owns c_1$, and $f^{\circ n}$ maps
$G_1,\ldots, G_{q-2}$ homeomorphically onto their images
in $G_2,\ldots, G_{q-1}$.
\end{lem}

As shown in \cite[Lemma 3.1]{BS1}, this characteristic point lies on
the arc $[0, c_1]$ and if $1 \IntAdr S_1 \IntAdr S_2 \IntAdr \dots$
is the internal address of $(T,f)$, then for every entry
$S_i$, there is a characteristic point $p_i$.

Recall that the map $\nu:\Circle \to \Sym$ assigns the
kneading sequence to an external parameter angle. 
In order to investigate how Hausdorff dimension behaves under $\nu^{-1}$,
we must determine, for an $n$-cylinder $C$, the number of components of
$\nu^{-1}(C)$
and their minimal length. This relies on the number of different ways
a Hubbard tree with an $m$-periodic critical point can be embedded
in the plane, because this equals
the number of components of $\nu^{-1}(C)$ for $m$-cylinders $C$.
Let $\varphi(q)$ be the Euler function counting the
integers $1 \leq i < q$ that are coprime to $q$; it gives the number of
transitive maps on $q$ points preserving circular order. 

\begin{lem}{(Embedding of the Hubbard Tree).}
\label{LemEmbedding} \lineclear
A Hubbard tree $(T,f)$ can be embedded into the plane so that $f$
respects the cyclic order of the local arms at all branch points if
and only if $(T,f)$ has no {\evil} orbits.
If $q_1,q_2,\ldots$ are the number of arms of the different
characteristic branch points (all of them {\tame}),
then there are $\prod \varphi(q_i)$
different ways to embed $T$ into the plane such that $f$ extends to
a two-fold branched covering.
\end{lem}

\begin{proof}
If $(T,f)$ has an embedding into the plane so that $f$ respects the
cyclic order of local arms at all branch points, then clearly there 
can be no {\evil} orbit
(this uses the fact that no periodic orbit of branch points contains 
a critical point).

Conversely, suppose that $(T,f)$ has no {\evil} orbits, so all local
arms at every periodic branch point are permuted transitively.
First we embed the arc $[0,c_1]$ into the plane, for example on a
straight line. Every cycle of branch points has at least its
characteristic point $p_1$ on the arc $[0,c_1]$, and it does not
contain the critical point. Suppose $p_1$ has $q$ arms. Take $s \in
\{ 1, \dots, q-1\}$ coprime to $q$
and embed the local arms at $p_1$ in such a way that
the return map $f^{\circ n}$ moves each arc over by $s$ arms
in counterclockwise direction.
This gives a single cycle for every $s < q$ coprime to $q$.
There are $\varphi(q)$ choices to do this and these choices can be made
for all characteristic branch points independently.

A point $x \in T$ is called \emph{marked} if it is a branch point or 
point on the critical orbit.
We say that two marked points $x,y$ are \emph{adjacent} if $(x,y)$
contains no further marked point. If a branch point $x$ is already
embedded together with all its local arms, and $y$ is an adjacent marked
point on $T$ which is not yet embedded but $f(y)$ is, then draw a
line segment representing $[x,y]$ into the plane, starting at $x$
and disjoint from the tree drawn so far. This is possible uniquely
up to homotopy. Embed the local arms at $y$ so that $f\colon y\to
f(y)$ respects the cyclic order of the local arms at $y$; this is
possible because $y$ is not the critical point of $f$.

Applying the previous step finitely many times, the entire tree $T$
can be embedded. It remains to check that for every characteristic
branch point $p_1$ of period $m$, say, the map $f\colon p_1\to
f(p_1)=:p_2$ respects the cyclic order of the local arms. By
construction, the forward orbit of $p_2$ up to its characteristic
point $p_1$ is embedded before embedding $p_2$, and
$f^{\circ(m-1)}\colon p_2\to p_1$ respects the cyclic order of the
embedding. If the orbit of $p_1$ is {\tame}, the cyclic order induced
by $f\colon p_1\to p_2$ (from the abstract tree) is the same as the 
one induced by $f^{\circ
(m-1)}\colon p_2\to p_1$ used in the construction (already embedded 
in the plane), and the
embedding is indeed possible.
Recalling that $q_1, q_2, \dots$ are the number of arms of the characteristic
branch points in $T$, we see that there are altogether $\prod \varphi(q_i)$
different ways to embed $T$.
\end{proof}

\begin{lem}{(Upper Bound for Number of Embeddings).}
\label{LemNumberEmbed2} \lineclear
A Hubbard tree in which the critical orbit is periodic with period
$n$ has less than $n$ embeddings into the plane that respect the
circular order of the local arms at every branch point.
\end{lem}
\begin{proof}
Let $1 \IntAdr S_1\IntAdr \dots \IntAdr S_k$ be the internal
address of the tree (cf.\ Definition~\ref{DefRho}), with $S_k=n$. We 
may suppose that all branch
points are {\tame} (or there would be no embedding at all).
By Lemma~\ref{LemHomeo}, the periods of all branch points appear
on the internal address.
Let $p_0, \dots , p_{k-1}$ be the {\tame} characteristic periodic
points of periods $S_0, \dots , S_{k-1}$. Let their numbers
of arms be $q_0, \dots, q_{k-1}$; according to
\cite[Proposition 4.19]{BS1} they satisfy
\[
S_{i+1} = \left\{ \begin{array}{ll}
(q_i-1)S_i + r_i &\mbox{ if } S_i \in \orb_{\rho}(r_i), \\
(q_i-2)S_i + r_i &\mbox{ if } S_i \notin \orb_{\rho}(r_i), \
\end{array} \right.
\]
where the $r_i$ are uniquely defined by the condition $1\leq
r_i\leq S_i$.

Since only branch points contribute to the number of embeddings,
let us write $i(0), i(1), \dots, i(l)$ for the indices of $p_i$
that are branch points. Obviously $k > i(l)$.

By Lemma~\ref{LemEmbedding}, there are precisely
$a := \prod_{s=0}^l \varphi(q_{i(s)})$ dynamically viable embeddings
of the Hubbard tree into the plane.
Clearly $a \leq \prod_{s=0}^l (q_{i(s)}-1)$.
We will show that $a < S_k$.
We call $\z_j$ be a {\em closest precritical point} of $\step(\z_j) = j$
if $f^{\circ j}(\z) = c_1$.
and the arc $[c_1, \z_j]$ contains no precritical point 
of $\step < j$.
The arc $[p_{i(t)}, c_1]$ contains the closest precritical point
$\z_{S_{i(t)+1}}$, and $f^{\circ (q_{i(t)}-2) S_{i(t)}}$
maps it to a precritical point $\z_t$ of
$\step(\z_t) = S_{i(t)+1} - (q_{i(t)}-2) S_{i(t)}$.
Lemma~\ref{LemHomeo} implies that the arm $G_1$ of $p_{i(t)}$
containing $c_1$
homeomorphically survives $f^{\circ (q_{i(t)}-2) S_{i(t)} }$
and $\z_t$ lies in a different arm of $p_{i(t)}$ as the critical point.
However, $\z_t$ and $\z_{S_{i(t-1)+1}}$ lie in the same global arm
of $p_{i(t-1)}$, which homeomorphically survives another
$(q_{i(t-1)}-2)S_{i(t-1)}$ iterates.
Inductively repeating this argument gives
\[
S_{i(t)+1} > (q_{i(t)}-2) S_{i(t)} + (q_{i(t-1)}-2) S_{i(t-1)} + \dots
+ (q_{i(0)}-2) S_{i(0)}.
\]
Choose $u_1 = S_{i(1)}$, $u_0 = u_1/(q_{i(0)}-2)$ and
\[
u_{t+1} := (q_{i(t)}-2) u_t + (q_{i(t-1)}-2) u_{t-1} + \dots
+ (q_{i(0)}-2) u_0
\,\,.
\]
Then by induction $u_{t+1} = (q_{i(t)}-1) u_t$, and therefore
$$
u_{t+1} = (q_{i(t)}-1) u_t = u_1 \prod_{s = 1}^t (q_{i(s)}-1).
$$
Hence $S_k \geq S_{i(l)+1} > u_{l+1} = S_{i(1)} \prod_{s=1}^l (q_{i(s)}-1)$.
It is easily checked that $S_{i(1)} \geq q_{i(0)} - 1$.
Therefore $S_k > a$ as asserted.
\end{proof}

\begin{proof}[Proof of Theorem~\ref{ThmHdim_mandel2}]
We know the dimension bounds in terms of kneading sequences $\nu$, which are
proved by means of counting $n$-cylinders.
Here we need to make the transition from parameter angle $\theta$
to $\nu(\theta)$. This involves counting how many arcs $A \subset \Circle$
map into the same cylinder set $C$ under $\nu$, which is related to how many 
ways there are to embed Hubbard trees into the plane.

For every $\theta$ with $\nu(\theta) \in C$, the Hubbard tree (whether finite 
or infinite) contains a finite skeleton composed of the connected hulls
of the characteristic periodic points of period up to $n$, see \cite{BKS}.
The number of possible embeddings of this skeleton coincides with the number of
different arcs in $\nu^{-1}(C)$, and hence we need to understand these embeddings only for finite trees.

For the upper bound, we claim that for every $n$-cylinder $C \subset \Symg$,
$\nu^{-1}(C)$ consists of at most
$\frac12 n(n+1)$ arcs of length $\leq \frac1{2^n-1}$.
Indeed, if $\alpha$ is such that the $n$-th entry $\nu(\alpha)_n = \*$, say
that $2^{n-1}\alpha = \frac{m+\alpha}{2}$
for some $m \geq 1$, then for $\alpha' = \alpha +  \frac1{2^n-1}$
we have $2^{n-1} \alpha' =  \frac{m+1+\alpha'}{2}$.
Therefore every component of $\nu^{-1}(C)$ must be contained
in an arc $(\alpha, \alpha +  \frac1{2^n-1})$ for some $\alpha \in
\Circle$. This shows that $\nu^{-1}$ is Lipschitz on each branch.

Let $T$ be a Hubbard tree with a periodic critical point;
say the period is $m = S_k$.
The external angles of $T$
depend on the specific embedding of $T$ in the plane. According to
Lemma~\ref{LemNumberEmbed2}
there are at most $m$ different embeddings.
Each embedding of $T$ (with biaccessible critical value)
comes with at least two external angles.
We can exclude the Hubbard trees with more than two external angles
at the critical value,
because these correspond to strictly preperiodic critical points
and this constitutes a countable set.
Hence there are at most $2m$ external angles realizing
the kneading sequence $\ovl{\nu_1\dots\nu_{m-1}\*}$.
Each arc in $\nu^{-1}(C_{\nu_1\dots\nu_n})$ has two boundary points
having kneading sequences $\ovl{\nu_1 \dots \nu_{m-1}\*}$
for some $m \leq n$.
Therefore the total number of arcs is bounded by
$\sum_{m=1}^n m = \frac12n(n+1)$. This proves the claim.
Now use Proposition~\ref{PropHdimBiacSymbolics} and 
Lemma~\ref{LemHdimPreserved} to finish the proof of the upper bound.

For the lower bound, take $M = N$ or $S_{\kappa+1}$ according to whether
the lower bound in Proposition~\ref{PropHdimBiacSymbolics} is 
obtained from $L_1(N)$ or $L_2(S_\kappa)$.
Let $k > M$ and take an $n$-cylinder set
$C = C_{e_1\dots e_n}$ intersecting $\Symg_{N,\kappa,k}$.
Without loss of generality we can choose $C$ so that
$n \in \orb_{\rho}(i)$ for each $i \leq M$, and that no $\nu \in C$ gives
rise to an {\evil} period $m$ with $\rho(m) \leq n$.

Using Proposition~\ref{PropHdimBiacSymbolics},
we can find a subset $B \subset C$ of Hausdorff dimension
$\delta \geq \max\{ L_1(N), L_2(S_\kappa) \}$.
Moreover, for all $\nu \in B$, $r_i := \rho(i) - i \leq M$ for all $i > n$.
Therefore $n \in \orb_{\rho}(r_i)$, so it follows
that every $\nu \in B$ corresponds to an admissible Hubbard tree
$T$, whose periodic branch points have period $\leq M$.
By Lemma~\ref{LemEmbedding},
$T$ has a bounded number of embeddings, hence the map
$\nu:\nu^{-1}(B) \to B$ is bounded-to-one.

A second property of $B \subset C$ is that
if $\tilde C = C_{e_1 \dots e_j}$
is any subcylinder intersecting $B$, then all four
subcylinders $C_{e_1 \dots e_je_{j+1}e_{j+2}}$ satisfy
Admissibility Condition~\ref{DefAdmissCond}.
Therefore the single arc component $A \subset \nu^{-1}(\tilde C)$ is
divided into four pieces by points of the form
$\frac{i}{2^{j+1}-1}$ or $\frac{i}{2^{j+2}-1}$ (where $i$ is an integer),
and $|A| > \frac1{2^{j+2}-1}$. It follows that the map
$\nu$ restricted to $\nu^{-1}(B)$ is Lipschitz (with Lipschitz constant $\leq 4$)
on each of its branches.
Therefore, we can use the second part of Lemma~\ref{LemHdimPreserved},
say with polynomial $Q(n) \equiv 4$, to conclude that
the set of biaccessible external angles
contains a Cantor set $\nu^{-1}(B)$ of Hausdorff dimension $\delta$.
\end{proof}

The following corollary deals with parameter angles whose rays land on 
hyperbolic components. Note that all these rays land indeed: if $\nu(\theta)$ is periodic of 
period $n$ but $\theta$ is irrational, then $\theta$ has a finite 
internal address and the parameter ray is contained in the wake of a 
hyperbolic component of period $n$, but not in any of its subwakes, and 
every boundary point of a hyperbolic component has trivial fiber 
\cite[Corollary 5.1]{Fibers2}.

\begin{cor}{(Hausdorff Dimension of Periodic Parameter Angles).}
\label{CorHausdorffParaAngles}\lineclear
The set of parameter angles $\theta\in \Circle$ such that $\nu(\theta)$
is periodic has zero Hausdorff dimension.
\end{cor}

\begin{proof} Let $\nu$ be a periodic kneading sequence. If $\nu$ is $\*$-periodic,
then there are at most finitely many $\theta\in \Circle$ such that
$\nu(\theta) = \nu$. Otherwise, there can be a Cantor set of such angles,
but the first half of the proof of Theorem~\ref{ThmHdim_mandel2}
shows that this Cantor set has Hausdorff dimension zero.
Since there are countably many periodic kneading sequences, the result
follows.
\end{proof}

\section{Biaccessibility and Renormalization}
\label{SecRenorm}

It remains to prove Propositions~\ref{PropDimRenorm} and \ref{PropDimBiacRenorm} on the biaccessibility dimension 
under the condition of renormalizability or being associated to the main molecule. 

\begin{proof}[Proof of Proposition~\ref{PropDimRenorm}]
If a kneading sequence is simple $M$-renormalizable
(regardless of whether it is admissible or not), then the associated
internal address contains after entry $M$ only entries which
are divisible by $M$, see \eqref{eq:nurenorm}. If the kneading sequence is divided into
blocks of length $M$, then every block can differ from the first one
only at the last position. Since $\nu_1 = \1$ and $\nu_M$ must be such that $M$
occurs in the internal address, there are at most $2^{M-2}\cdot
2^{i-1}$ possibilities for the first $iM$ entries of
such kneading sequences. Hence the set of
$M$-renormalizable kneading sequences has Hausdorff dimension
at most $(\log 2)/\log(2^M)=1/M$ by Lemma~\ref{LemBasic_hdim}. 
Infinitely renormalizable kneading sequences are simple $M$-renormalizable 
for arbitrarily large $M$, so their Hausdorff dimension is $0$.

Given a parameter angle $\varphi$ with kneading sequence $\nu(\varphi)$, define the interval
\begin{eqnarray*}
J_n(\varphi) = \{ \varphi' \in \Circle &:& \nu(\theta)_i = \nu(\varphi)_i 
\text{ for all } 1 \leq i < n \text{ and } \\
&&  \theta \text{ on the shorter arc } (\varphi, \varphi')\}.
\end{eqnarray*}
Then $|2^i \theta - \2^i \theta'| < 1$ for all $\theta,\theta' \in J_n(\varphi)$ 
and $1 \leq i < n$, so $\diam(J_n(\varphi)) \leq 2^{1-n}$.

If $\varphi$ is $M$-renormalizable, then, as above, there at most $2^{M-2} 2^{i-1}$ ways 
to select the first $iM$ digits.
By Lemma~\ref{LemNumberEmbed2}, there are at most $iM$ ways in which the 
corresponding Hubbard tree can be embedded in the plane.
Hence, the set of $M$-renormalizable parameter angles is covered by
a collection of at most $iM 2^{M-3+i}$ intervals of length $2^{1-iM}$.
Therefore the Hausdorff dimension of the set of
$M$-renormalizable parameter angles is bounded by
$\liminf_i \frac{ \log(iM 2^{M-3+i}) }{ \log 2^{1-iM} } = \frac1M$.
\hide{
If a parameter angle $\varphi$ (with kneading sequence $\nu$)
is $M$-renormalizable, then, for each $i \geq 1$,
$\varphi$ is separated from $0$ by a ray-pair
$\pair{\theta}{\theta'}$ with common kneading sequence 
$\ovl{ \nu_1\dots\nu_{Mi-1}\* }$, and $|\theta-\theta'| \le 1/(2^{Mi}+1)$.
As above, there at most $2^{M-2} 2^{i-1}$ ways to select the first $iM$ digits.
By Lemma~\ref{LemNumberEmbed2}, there are at most $iM$ ways in which the 
corresponding Hubbard tree can be embedded in the plane.
Hence, the set of $M$-renormalizable parameter angles is covered by
a collection of at most $iM 2^{M-3+i}$ intervals of length $1/(2^{Mi}+1)$.
Therefore the Hausdorff dimension of
$M$-renormalizable parameter angles is bounded by
$\liminf_i \frac{ \log(iM 2^{M-2+i}) }{ \log(2^{iM}+1) } = \frac1M$.
}
For infinitely renormalizable angles, this holds for arbitrary
large $M$, so the Hausdorff dimension of this set is $0$.
\end{proof}

\begin{proof}[Proof of Proposition~\ref{PropDimBiacRenorm}]
If $\theta$ is associated to the main molecule of $\M$ and infinitely renormalizable, then $\kappa = \infty$ and $\dim_H(\Biac_\theta) \le U_2(S_\kappa) = 0$
by Theorem~\ref{ThmHdimBiaccAngles}.
If $\theta$ is associated to the main molecule and finitely simple renormalizable, then in particular, it has a finite internal address, and 
 there are only countably many ray-pairs in this case.

The other implication follows immediately from
Theorem~\ref{ThmHdimBiaccAngles}
because if $\theta$ is not associated to the main molecule of $\M$, 
then $S_\kappa < S_{\kappa+1} < \infty$ and we have the lower
bound $L_2(S_\kappa) = 1/S_{\kappa+1} > 0$.
\end{proof}


\begin{thebibliography}{BGMY}
\normalsize 

\bibitem[BKS]{BKS} Henk Bruin, Alexandra Kaffl, Dierk Schleicher,
{\em  Existence of quadratic Hubbard trees,}
Fund.\ Math.\ {\bf 202} (2009), 251--279.

\bibitem[BrS]{BS1} Henk Bruin, Dierk Schleicher,
{\em Admissibility of kneading sequences and structure of Hubbard 
trees for quadratic polynomials,}
J.\ London.\ Math.\ Soc.\ {\bf 8} (2009),  502--522.

\bibitem[BuS]{BulSen} Shaun~Bullett, Pierrette~Sentenac,
{\em Ordered orbits of the shift, square roots,
and the devil's staircase,}
Math. Proc. Cambridge Philos. Soc. {\bf 115} (1994),
451--481.

\bibitem[CT]{CT} Carlo Carminati, Giulio Tiozzo,
{\em The bifurcation locus of numbers of bounded type,}
Preprint 2012, arXiv:1111.2554.

\bibitem[C]{Chung} Kai Lai Chung,
{\em A course in probability theory,}
Harcourt, Brace \& World, Inc., New York (1968).

\bibitem[D]{DouadyCompacts} Adrien Douady, 
\emph{Descriptions of compact sets in $\C$}, in:
Topological Methods in Modern Mathematics, Publish or Perish
(1993), 429--465.

\bibitem[DH]{Orsay} Adrien~Douady, John~Hubbard,
{\em \'Etudes dynamique des polyn\^omes complexes
I \& II},
Publ. Math. Orsay. (1984-85) {\em(The Orsay notes)}.

\bibitem[DS]{DS} Dzmitry Dudko, Dierk Schleicher, 
\emph{Core entropy of quadratic polynomials}. With an appendix by Wolf Jung.  
Preprint 2014, arXiv:1412.8760.

\bibitem[HY]{HubbardYoccoz}
John Hubbard, \emph{Local connectivity of bifurcation loci: three theorems of Jean-Christophe Yoccoz}. In: Topological Methods in Modern Mathematics. Publish or Perish, Houston, TX 1993, {375--378} and {467--511}.

\bibitem[J]{jung} Wolf Jung, 
{\em Core entropy and biaccessibility of quadratic polynomials,}
Preprint 2014,  arXiv:1401.4792 

\bibitem[KL1]{KahnLyubich}
Jeremy Kahn, Mikhail Lyubich, 
\emph{A priori bounds for some infinitely renormalizable quadratics: II. Decorations}.  
Annales Scientifiques de l'Ecole Normale Sup\'erieure.
 {\bf 41} (2008), 57--84. 

\bibitem[KL2]{KahnLyubich2}
Jeremy Kahn, Mikhail Lyubich, 
\emph{A priori bounds for some infinitely renormalizable quadratics: III. Molecules}.  
In: \emph{Complex dynamics: families and friends} (ed.\ D.\ Schleicher), 2009. 
 
\bibitem[LS]{IntAdr} Eike~Lau, Dierk~Schleicher,
{\em Internal addresses in the Mandelbrot set and irreducibility
of polynomials,}
Stony Brook Preprint {\bf \#19} (1994).


\bibitem[La]{Lavaurs} Pierre Lavaurs,
{\em Une description combinatoire de l'involution d\'efinie par
$M$ sur les rationnels \`a d\'enominateur impair},
C. R. Acad. Sci. Paris, S\'erie I Math. {\bf 303} (1986),
143--146.

\bibitem[Ly1]{LyubichQuadratics}
Mikhail Lyubich, \emph{Dynamics of quadratic polynomials, I-II}, Acta Math., \textbf{178} (1997), 185--297.


\bibitem[Ly2]{Lyubich}  Mikhail~Lyubich,
{\em How big is the set of infinitely renormalizable quadratics?}
In ``Voronezh Winter Mathematics School,'' Amer.\ Math.\ Soc.\ 
Transl.\ Ser.\ 2, {\bf 184}, Amer.\ Math.\ Soc., Providence, RI, 
(1998), 131--143.

\bibitem[Man]{Manning} Anthony~Manning,
{\em Logarithmic capacity and renormalizability for landing on
the Mandelbrot set,}
Bull. London Math. Soc. {\bf 28} (1996), 521--526.

\bibitem[Mat]{Mattila} Pertti Mattila,
{\em Geometry of sets and measures in Euclidean spaces,}
Cambridge University Press, (1995).

\bibitem[Mi2]{MiOrbits} John Milnor, 
\emph{Periodic orbits, external rays, and the Mandelbrot
set: an expository account}. Ast\'erisque {\bf 261} (2000),
277--333.

\bibitem[MS]{MeerkampSchleicher}  Philipp Meerkamp, Dierk Schleicher,
{\em Hausdorff dimension and biaccessibility for polynomial Julia sets,}
Proc.\ Amer.\ Math.\ Soc.\ {\bf 141} (2013), 533--542.

\bibitem[P]{Pe} Chris~Penrose,
{\em On quotients of shifts associated with dendrite Julia sets
of quadratic polynomials,}
Ph.D.\ Thesis, University of Coventry, (1994).

\bibitem[Pet]{Petersen}
Carsten Lunde Petersen, 
\emph{Local connectivity of some Julia sets containing a circle with an irrational rotation}, Acta. Math., \textbf{177} (1996), 163--224. 

\bibitem[R]{rees} Mary~Rees,
{\em  A partial description of the parameter space of rational maps of
degree two: Part 1,}
Acta. Math. {\bf 168} (1992), 11--87
(See also: Realization of matings of polynomials as rational maps of degree
two, Preprint 1986).

\bibitem[Sch1]{IntAddr} Dierk Schleicher, \emph{Internal Addresses in 
the Mandelbrot Set and Galois Groups of Polynomials}. Preprint 2012, arXiv:9411238. 
Arnold Mathematical Journal, 
appeared online Aug.\ 2016  DOI 10.1007.540598-016-0042-x

\bibitem[Sch2]{Fibers1} Dierk Schleicher, 
\emph{On fibers and local connectivity of compact sets in $\C$}. 
Stony Brook preprint \textbf{12} (1998). 
arXiv:math/9902154

\bibitem[Sch3]{Fibers3} Dierk Schleicher, 
\emph{On fibers and  renormalization of Julia sets and Multibrot sets}. 
Stony Brook Preprint \textbf{13b} (1998). 
arXiv:math/9902156

\bibitem[Sch4]{ExtRayMandel} Dierk Schleicher,
\emph{Rational external rays of the Mandelbrot set},
Asterisque {\bf 261} (2000), 405--443.

\bibitem[Sch5]{Fibers2} Dierk~Schleicher,
{\em On fibers and local connectivity of Mandelbrot and Multibrot
sets}. In: M.~Lapidus, M.~van Frankenhuysen (eds): {\em
Fractal Geometry and Applications: A Jubilee of Beno\^{\i}t Mandelbrot}.
Proceedings of Symposia in Pure Mathematics {\bf 72},
American Mathematical Society (2004), 477--507.

\bibitem[SZ]{SchleicherZakeri}
Dierk Schleicher, Saeed Zakeri,
\emph{On biaccessible points in the Julia set of a Cremer quadratic 
polynomial}. Proc. Amer. Math. Soc \textbf{128} 3 (1999), 933--937.

\bibitem[Sh]{shishikura} Mitsuhiro~Shishikura,
{\em On a theorem of M. Rees for matings of polynomials}, in:
Tan Lei (ed.), {\em The Mandelbrot set, theme and variations},
Cambridge University Press {\bf 274} (2000), 289--305.

\bibitem[Si1]{sirvent1} V\'{\i}ctor Sirvent,
{\em Space-filling curves and geodesic laminations,}
Geom.\ Dedicata, {\bf 135} (2008), 1--14.

\bibitem[Si2]{sirvent2} V\'{\i}ctor Sirvent,
{\em Space-filling curves and geodesic laminations II, symmetries,}
Monats.\ Math.\ {\bf 166} (2012), 543--558.

\bibitem[Sm]{smirnov1} Stanislav~Smirnov,
{\em On the support of dynamical laminations and biaccessible points 
in Julia set,}
Colloq. Math. {\bf 87} (2001), 287--295.

\bibitem[Ta1]{tanlei1} Tan Lei,
{\em Similarity between the Mandelbrot set and Julia sets,}
Commun. Math. Phys. {\bf 134} (1990), 587--617.

\bibitem[Ta2]{tanlei2} Tan Lei,
{\em Matings of quadratic polynomials,}
Ergod. Th. \& Dynam. Sys. {\bf 12} (1992), 589--620.

\bibitem[Ta3]{TanLeiEntropy} Tan Lei, \emph{On W.\ Thurston's core-entropy theory}. Presentation, given in Toulouse (January 2014) and elsewhere.
 
\bibitem[Thn]{thunberg} Hans~Thunberg,
{\em A recycled characterization of kneading sequences,}
Internat. J. Bifur. Chaos Appl. Sci. Engrg.  {\bf 9}  (1999), 1883--1887.

\bibitem[Th1]{ThurstonLaminations} William P.\ Thurston, \emph{On the geometry and dynamics of iterated rational maps}. In: Complex dynamics, families and friends (Dierk Schleicher, ed.), AK Peters, Wellesley, MA (ISBN 978-1-56881-450-6), pp. 3--109 (2009).

\bibitem[Th2]{thurston} William Thurston, 
{\em Entropy in dimension one,} 
Frontiers in complex dynamics, 339--384,  Edited by S.\ Koch, Princeton Math. Ser., {\bf 51}, Princeton Univ. Press, Princeton, NJ, 2014.

\bibitem[Ti1]{tiozzo} Giulio Tiozzo,
{\em Topological entropy of quadratic polynomials and dimension of sections of the Mandelbrot set,}
 Adv.\ Math.\ {\bf 273} (2015), 651--715.
   
\bibitem[Ti2]{tiozzo2} Giulio Tiozzo,
\emph{Continuity of core entropy of quadratic polynomials}. 
Inventiones Mathematicae, {\bf 203} (2016), no.\ 3, 891--921.

\bibitem[Za]{zakeri3} Saeed~Zakeri,
{\em Biaccessibility in quadratic Julia sets,}
Ergod. Th. \& Dynam. Sys. {\bf 20} (2000), 1859--1883.

\bibitem[Zd]{zdunik} Anna~Zdunik,
{\em On biaccessible points in Julia sets of polynomials,}
Fund. Math. {\bf 163} (2000), 277--286.

\end{thebibliography}
\end{document}